%% file: geomed.tex
\documentclass[review,onefignum,onetabnum]{siamonline190516}

\input{geomedShared}

\ifpdf
\hypersetup{
  pdftitle={The Geometric Median and Applications to Robust Mean Estimation},
  pdfauthor={S. Minsker, N. Strawn}
}
\fi

\usepackage{mathtools}
\usepackage{amsmath,amssymb}
\usepackage{subfig, float}
\usepackage{graphicx}
\usepackage{verbatim,epsfig,color,lscape}
\usepackage{epstopdf}
\usepackage{color}
\def\R{\mathbb{R}}
\def\E{\mathbb{E}}
\def\P{\mathbb{P}}

\def\l{\left}
\def\r{\right}
\newcommand{\be}[1]{\begin{equation*}#1\end{equation*}}
\newcommand{\ben}[1]{\begin{equation}#1\end{equation}}
\def\al#1{\begin{align*}{#1}\end{align*}}

\def\ml#1{\begin{multline*}{#1}\end{multline*}}
\def\mln#1{\begin{multline}{#1}\end{multline}}

\newcommand{\wt}{\widetilde}
\newcommand{\wh}{\widehat}

\newcommand{\m}{\mathcal}
\newcommand{\mb}{\mathbb}
\newcommand\argmin{\mathop{\mbox{argmin}}}

\newcommand{\tr}{\mbox{tr}\,}
\newcommand{\var}{\mbox{Var}}

\newcommand{\eps}{\varepsilon}

\newcommand{\vol}[1]{\mbox{vol}\left(#1\right)}
\newcommand{\med}[1]{\mbox{med}\left(#1\right)}

\newcommand{\pr}[1]{\mathbb{P}{\left(#1\right)}}
\newcommand{\dotp}[2]{\left\langle#1,#2\right\rangle}

\nolinenumbers

\begin{document}

\maketitle

\begin{abstract}
This paper is devoted to the statistical and numerical properties of the geometric median, and its applications to the problem of robust mean estimation via the median of means principle. Our main theoretical results include (a) an upper bound for the distance between the mean and the median for general absolutely continuous distributions in $\mathbb R^d$, and examples of specific classes of distributions for which these bounds do not depend on the ambient dimension $d$; (b) exponential deviation inequalities for the distance between the sample and the population versions of the geometric median, which again depend only on the trace-type quantities and not on the ambient dimension. As a corollary, we deduce improved bounds for the (geometric) median of means estimator that hold for large classes of heavy-tailed distributions. 
Finally, we address the error of numerical approximation, which is an important practical aspect of any statistical estimation procedure. We demonstrate that the objective function minimized by the geometric median satisfies a ``local quadratic growth'' condition that allows one to translate suboptimality bounds for the objective function to the corresponding bounds for the numerical approximation to the median itself, and propose a simple stopping rule applicable to any optimization method which yields explicit error guarantees. 
We conclude with the numerical experiments including the application to estimation of mean values of log-returns for S\&P 500 data.
\end{abstract}



\section{Introduction}
\label{sec:intro}

The geometric median, also referred to as the spatial median and the $L_1$ median, is one of the oldest and most popular robust estimators of location. Its roots go back to the Fermat, Toricelli and Weber \cite{weber1929uber} under the name of ``Fermat-Weber'' point, and to the work \cite{haldane1948note} under the name of Haldane's median; other notable early references include the paper by Gini and Galvani \cite{gini1929di}. The geometric median is an element of a more general family of spatial quantiles that was introduced and studied in detail by Koltchinskii and Chaudhuri \cite{koltchinskii1994spatial,Koltchinskii1997M-estimation-co00,chaudhuri1996geometric}: in particular, existence, uniqueness, and the asymptotic properties of spatial quantiles are well-understood. Extensions of the geometric median to the general Banach spaces were analyzed by Kempreman \cite{kemperman1987median} and, more recently, by Romon \cite{romon2022statistical}. Deep connections between the probability distributions and the corresponding spatial quantiles have been investigated by Konen \cite{konen2022recovering}. 

Renewed interest in the properties of the geometric median was sparked with the re-introduction of the so-called ``median of means'' (MOM) estimator into high-dimensional statistics and machine learning literature. Originally appearing in the works of Nemirovski and Yudin \cite{Nemirovski1983Problem-complex00,jerrum1986random,alon1996space} in a different context, the MOM estimator was shown to be a powerful tool for the analysis of corrupted and heavy-tailed data by Lerasle and Oliveira \cite{lerasle2011robust}. The work by Hsu and Sabato \cite{hsu2016loss} demonstrated multiple novel applications of the original estimator by Nemirovski and Yudin in general metric spaces, while Minsker \cite{minsker2015geometric} introduced a version of the median of means principle based on the geometric median. On a high level, the median of means estimator can be viewed as a ``majority vote'' among several independent estimators of the mean. Its popularity can be attributed to the fact that it is widely applicable, efficiently computable even in high dimensions, requires minimal tuning, and admits strong theoretical guarantees in many circumstances. 
However, as was pointed out by several authors, for instance by Lugosi and Mendelson \cite{lugosi2017sub}, the geometric median of means estimator fails to attain optimal deviation bounds for the fundamental problem of multivariate mean estimation. Specifically, let $Y_1,\ldots,Y_N$ be i.i.d. copies of a random vector $Y\in \mb R^d$ with mean $\mb EY=\mu$ and covariance $\mb E(Y-\mu)(Y-\mu)^T = \Sigma_Y$. Then, as shown by Minsker \cite{minsker2015geometric}, for any $1\leq t\leq N/2$, there exists a version $\wh\mu_N=\wh\mu_N(Y_1,\ldots,Y_N;t)$ of the geometric MOM estimator, formally defined in display \eqref{eq:geom-med} below, such that 
\ben{
\label{eq:mom-basic}
\l\| \wh\mu_N - \mu\r\|\leq C\sqrt{\frac{\tr(\Sigma_Y)t}{N}}
}
with probability at least $1-e^{-t}$; here, $C>0$ is an absolute constant, $\|\cdot\|$ stands for the Euclidean norm of a vector and the spectral norm of a matrix, and $\tr(\cdot)$ denotes the trace of the operator. At the same time, a sub-Gaussian estimator $\wt \mu_N$ should satisfy an inequality akin to the sample mean of a Gaussian distribution, namely,
\ben{
\label{eq:sg}
\l\| \wt\mu_N - \mu\r\|\leq C\l( \sqrt{\frac{\tr(\Sigma_Y)}{N}} +  \sqrt{\|\Sigma_Y\|}\sqrt{\frac{t}{N}}\r)
}
with probability at least $1-e^{-t}$, where $C>0$ is an absolute constant; the advantage of the latter inequality over \eqref{eq:mom-basic} is the fact that the deviation parameter $t$ and the dimension-dependent quantity $\tr(\Sigma_Y)$ appear in separate additive terms. It immediately implies that the radii of the confidence balls for the true mean $\mu$ derived from the inequality \eqref{eq:sg} are much smaller compared to their counterpart obtained from \eqref{eq:mom-basic}. Lugosi and Mendelson \cite{lugosi2017sub} proposed an alternative to the standard median of means principle based on the notion of tournaments and showed that the resulting estimator achieves the desired sub-Gaussian deviation guarantees for distributions possessing only the finite second moment. Many improvements, extensions and refinements of sub-Gaussian estimators have been suggesteed in the mathematical statistics and theoretical computer science literature since: we refer the reader to the excellent surveys by Lugosi and Mendelson \cite{lugosi2019mean} and Diakonikolas and Kane \cite{roughgarden2021beyond}. While the original estimator by Lugosi and Mendelson \cite{lugosi2017sub} is difficult to compute, several closely related numerically feasible alternatives have been proposed by Hopkins \cite{hopkins2020mean}, Cherapanamjeri \cite{cherapanamjeri2019fast}, Depersin and Lecu\'{e} \cite{depersin2022robust}, Bateni et al. \cite{bateni2022nearly}, among others. However, to the best of our knowledge, none of these methods admit practical implementations comparable to the best algorithms for evaluating the geometric median, as those in the works by Cohen et al. \cite{cohen2016geometric}, Beck and Sabach \cite{beck2015weiszfeld} and Cardot et al. \cite{cardot2017online}. 
Due to the computational advantages offered by the geometric median of means, it has become a popular tool for designing robust versions of distributed optimization methods such as Federated Learning \cite{alistarh2018byzantine,bouhata2022byzantine,chen2017distributed,pillutla2022robust}. Therefore, improved guarantees for the geometric MOM estimator have immediate implications for a variety of algorithms that use MOM principle as a subroutine. 

\subsection{Statistical error bounds}

In this paper, we revisit the original geometric median of means construction and show that the inequality \eqref{eq:mom-basic} can be improved for large classes of absolutely continuous, heavy-tailed distributions with sufficiently large effective rank $r(\Sigma_Y):=\frac{\tr(\Sigma_Y)}{\|\Sigma_Y\|}$ (indeed, if $r(\Sigma_Y)$ is bounded by a constant, then \eqref{eq:mom-basic} readily provides sub-Gaussian guarantees). Specifically, we show that $\wh\mu_N$ satisfies the bound of sub-exponential type: for all $t\lesssim \sqrt N$ (where $\lesssim$ denotes the inequality up to an absolute multiplicative constant), there exists a version of the MOM estimator $\wh\mu_N$ such that
\ben{
\label{eq:se}
\l\| \wh\mu_N - \mu\r\|\leq C\l( \sqrt{\frac{\tr(\Sigma_Y)}{N}} +  \sqrt{\|\Sigma_Y\|}\frac{t}{\sqrt N}\r)
}
with probability at least $1-e^{-t}$. While this bound succeeds in separating the trace and the confidence parameter $t$ into different additive terms as in \eqref{eq:sg}, thus making a significant improvement over the previously known bound \eqref{eq:mom-basic}, it fails to achieve optimal sub-Gaussian behavior. A remaining open question is whether there exists an easily described class of heavy-tailed distributions for which the geometric median admits truly sub-Gaussian deviation bounds. 

While the proof of the inequality \eqref{eq:mom-basic} is based on a simple ``majority vote-type'' argument, the present analysis leading to \eqref{eq:se} blends accurate estimates for the bias and the stochastic error of the geometric median of means. The upper bound for the bias (Theorem \ref{th:bias} and section \ref{section:moments}) is shown to be controlled by ratios of the negative moments of the norm that in turn depend on the ``small ball'' probability estimates. Control of the stochastic error relies on the deviation bounds for the geometric median (Theorem \ref{th:stochastic}) that, to the best of our knowledge, are new. In particular, our bounds depend only on trace-type quantities and not on the dimension of the ambient space, and yield sub-Gaussian type guarantees for a wide range of confidence levels.
 
\subsection{Numerical error bounds} 

Recall that the geometric median associated with the distribution $P_Y$ of a random vector $Y\in \mb R^d$ is defined as  
\[
m(P_Y) := \argmin\limits_{z\in \mb R^d} \mb E\l( \|z-Y\| - \|Y\|\r).
\]
Its empirical version based on an i.i.d. sample $Y_1,\ldots,Y_k$ is 
\ben{
\label{eq:geom-med}
\wh m = \med{Y_1,\ldots,Y_k}:=\argmin\limits_{z\in \mb R^d} \frac{1}{k}\sum_{j=1}^k \|z-Y_j\|.
}
In the sequel, we will frequently write $F(z) = F(z;Y_1,\ldots,Y_k)$ in place of $\frac{1}{k}\sum_{j=1}^k \|z-Y_j\|$. It is well-known that $m$ and $\wh m$ are well-defined and are unique unless $P_Y$ (or its empirical counterpart $\wh P_k = \frac{1}{k}\sum_{j=1}^k \delta_{Y_j}$) is supported on a straight line. 

Theoretical guarantees usually describe the performance of the ``ideal'' estimator $\wh m$ that is never known exactly. In practical applications, the ability to quantify numerical error of the algorithms used to approximate $\wh m$ gives one a litmus test for the overall performance of the estimator. However, most known results focus on suboptimality bounds for the objective function $F(z)$, while explicit error bounds for the approximation to the median itself are not available, to the best of our knowledge. 
More specifically, convex optimization theory usually asks questions about the computational complexity of finding a point $z_\eps$ such that $F(z_\eps)\leq F(\wh m) + \varepsilon$ for a given threshold $\varepsilon>0$. Existing results are fully quantitative, but for statistical applications we require the bounds for $\Vert z - \wh m\Vert$ instead. To bridge this gap, we prove (in Theorem \ref{thm:qgc}) a quadratic growth condition of the form
\[
F(z)-F(\wh m) \geq C_1 \frac{\Vert z - \wh m\Vert^2}{\Vert z - \wh m\Vert + C_2},
\]
where $C_1$ and $C_2$ are explicit functions of the data $Y_1,\ldots,Y_k$. This inequality immediately promotes any sub-optimality bound for the objective function to an error bound for approximating the median. As a corollary, we deduce a practical stopping criteria for any algorithm designed to find the geometric median, and propose a simple numerical procedure with fully explicit error bounds. 

The problem of minimizing $F(z)$ is a classical one, and has a long history. The most well-known numerical method is perhaps the celebrated \textit{Weiszfeld's algorithm} \cite{Weiszfeld1936Sur-un-probleme00}. Various improvements, refinements and accelerated versions of Weiszfeld's algorithm have been proposed and analyzed over the years. For example, results in this direction have been obtained in \cite{ostresh1978convergence,overton1983quadratically,karkkainen2005computation,vardi2000multivariate}, among others; an excellent review of the state of the art along with several new advances is given by Beck and Sabach \cite{beck2015weiszfeld}. The work by Cardot et al. \cite{cardot2013efficient} develops an online stochastic descent algorithm for minimizing $F(z)$ and provide its asymptotic convergence rate. The interior point method with the best-to-date convergence guarantees has been developed by Cohen et al. \cite{cohen2016geometric}; this work also provides a thorough comparison of existing alternatives. 

\subsection{Organization}

The rest of the content is organized as follows: in Section \ref{sec:main}, we introduce key notation and state our main results. Section \ref{sec:statistical} is devoted to the non-asymptotic analysis of the statistical properties of the geometric median, and culminates in the proof of Theorem \ref{th:main}. In section \ref{sec:numerical}, we discuss numerical algorithms for computing the geometric median that admit quantifiable error bounds, and prove a so-called quadratic growth condition. Section \ref{sec:exp} concludes the paper with the numerical experiments.

\section{Main results}
\label{sec:main}

Let us recall the definition of the median of means estimator based on a sample $Y_1,\ldots,Y_N$. Let $G_1\cup\ldots\cup G_k\subseteq \l\{ 1,\ldots,N\r\}$ be an arbitrary collection of $k\leq N/2$ disjoint subsets (``blocks'') of cardinality $n=\lfloor N/k\rfloor$ each, $\bar Y_j:=\frac{1}{|G_j|}\sum_{i\in G_j} Y_i$ and 
\ben{
\label{eq:geom-mom}
\wh \mu_N := \med{ \bar Y_1,\ldots,\bar Y_k}.
}
The main goal of this work is to understand when the random variable $\l\| \wh\mu_N - \mu\r\|$ admits good deviation bounds under minimal assumptions on the distribution of $Y$, and what is the typical computational complexity of approximating $\wh\mu_N$. We will now define the classes of distributions for which such ``good bounds'' can be established. 

Everywhere below, it will be assumed that the distribution of a random vector $Y$ is absolutely continuous with respect to the volume measure on a linear subspace of $\mb R^d$ (the linear span of the support of $P_Y$), and $M(Y)$ will stand for the sup-norm of the corresponding density $p_Y$. Similarly, if $X\in \mb R$ is a random variable with absolutely continuous distribution, $M(X)$ will denote the sup-norm of its density. 
The classes of distributions we are interested in are defined next.

\begin{enumerate}
\item \textbf{Linear transformations of the independent factors:} let $Y\in \mb R^d$ be given by a linear transformation $Y = AX$ where $X=(X_1,\ldots,X_k)\in \mb R^{k}$ is a centered random vector with independent coordinates such that $\Sigma_X$ is the identity matrix $I_k$ and $\max_{j=1,\ldots,k} M(X_j) =: M_0 <\infty$. Moreover, assume that $\max_{j=1,\ldots,k}\mb E|X_j - \mb EX_j|^q = K(q)<\infty$ for some $q>2$. The class of corresponding distributions $P_Y$ will be denoted $\m P_1:=\m P_1(M_0,K)$.  

\item \textbf{Distributions with well-conditioned covariance matrices:} let $Y\in \mb R^d$ be a random vector with support contained in a k-dimensional subspace $L$ such that its distribution is absolutely continuous with respect to the volume measure on $L$. Assume that \footnote{Here, we implicitly view $\Sigma_Y$ as an operator $\Sigma_Y: L\mapsto $L.}
\begin{enumerate}
\item $M^{1/k}\l( \Sigma^{-1/2}_Y Y\r)\leq M_0$;
\item $\frac{\tr(\Sigma_Y)}{ k \cdot\det^{1/k}(\Sigma_Y)}\leq R$;
\item For some $q>2$ and all unit vectors $u$, 
\begin{align*}
\label{eq:moment-equiv}
\mb E^{2/q} \l| \dotp{Y}{u} \r|^q \leq K\dotp{\Sigma_Y u}{u}.
\end{align*}
\end{enumerate}
The class of all such distributions will be denoted $\m P_2:=\m P_2(k,M_0,K,R)$.
\item \textbf{Signal plus noise:} let $Y = X + \xi\in \mb R^d$ where $P_X\in \m P_2(k,M_0,K,R)$, $\xi$ is independent from $X$ and is such that $\tr (\Sigma_\xi) \leq h\,\tr(\Sigma_X)$. This class of distributions is a natural generalization of $\m P_2(k,M_0,K,R)$ and will be denoted $\m P_3:=\m P_3(k,M_0,K,R,h)$. Distributions from the class $\m P_3$ can naturally be viewed as perturbations of the elements of the class $\m P_2$.
\end{enumerate}
The first main result of the paper is the following high-probability bound for the estimator $\wh \mu_N$.
\begin{theorem}
\label{th:main}
Assume that the distribution of $Y$ belongs to the class $\m P_j, \ j\in\{1,2,3\}$. Then for all $k_0\leq k\leq N/2$, the median of means estimator $\wh\mu_N$ defined in \eqref{eq:geom-mom} satisfies the inequality
\ben{
\label{eq:main-bound}
\|\wh \mu_N - \mu\|\leq C\l(\sqrt{\frac{\tr(\Sigma_Y)}{N}} + \sqrt{\|\Sigma_Y\|}\sqrt{\frac{k}{N}}\r)
}
with probability at least $1-e^{-\sqrt k}$, where $k_0$ and $C$ depend only on the parameters of the corresponding class $\m P_j$.
\end{theorem}

\begin{remark} 
\label{remark:key}
Let us discuss the main assumptions of the theorem.
\begin{enumerate}
\item Note that $\frac{\tr(\Sigma_Y)}{ k\cdot\det^{1/k}(\Sigma_Y)}$ is the ratio of the arithmetic and the geometric means of the eigenvalues $\lambda_1\geq\ldots\geq \lambda_k$ of $\Sigma_Y$: this quantity behaves well when the eigenvalues are of ``similar'' magnitude. For example, if $\lambda_j = \frac{C}{j^\alpha}$ for $\alpha<1$, then it is easy to check that 
\[
\frac{\sum_{j=1}^k\lambda_j}{k\l(\prod_{i=1}^k \lambda_i\r)^{1/k}} \leq C(\alpha).
\]
In fact, it is known \cite{gluskin2003note} that for most (with respect to the uniform distribution on a sphere) sequences, the ratio of arithmetic and geometric means is well-behaved.

\item Moment equivalence conditions similar to \eqref{eq:moment-equiv} are well known in the literature - for example, it has been employed in \cite{mendelson2020robust,lugosi2020multivariate,zhivotovskiy2021dimension,oliveira2016lower}, among others, in the contexts of robust estimation and random matrix theory. It is known to hold (Lemma 4.2 in \cite{mendelson2015learning}) for random vectors of the form $Y = AX$ where $X$ is either a vector with independent coordinates, or an unconditional vector with coordinates possessing finite moments of order $q$ (recall that a random vector has unconditional distribution when the distribution of  $(\eps_1X_1,\ldots,\eps_dX_d)$ is the same as the distribution of $X=(X_1,\ldots,X_d)$ for any sequence 
$\eps_1,\ldots,\eps_d\in \{\pm1\}^d$. Many elliptically symmetric distributions, for example multivariate Student's t-distribution, also satisfy \eqref{eq:moment-equiv} under appropriate restrictions on the number of degrees of freedom. 

Define the spatial sign covariance matrix via $D_Y:=\mb E\l[ \frac{(Y-m)}{\|Y-m\|} \frac{(Y-m)^T}{\|Y-m\|}\r]$, where $m=m(P_Y)$ is the geometric median of $Y$. The role of assumption \eqref{eq:moment-equiv} is in showing that $\Delta:=\|D_Y\| \leq \frac{C}{r(\Sigma_Y)}$. When \eqref{eq:moment-equiv} does not hold, inequality \eqref{eq:main-bound} is still valid with $\sqrt{\|\Sigma_Y\|}$ replaced by $\max\l(\sqrt{\|\Sigma_Y\|},\sqrt{\frac{\tr(\Sigma_Y)\Delta}{\sqrt k}} \r)$.
\end{enumerate}
\end{remark}
In the following sections, we develop the technical tools needed to prove Theorem \ref{th:main} and discuss the numerical methods used to approximate the estimator $\wh\mu_N$.
The proof of Theorem \ref{th:main} is based on the error decomposition 
\ben{
\label{eq:bound1}
\l\|\wh \mu_N - \mu\r\| \leq \l\| m_n - \mu \r\| + \l\|\wh\mu_N - m_n\r\|
}
where $m_n$ is the geometric median of the distribution $P^{(n)}$ of the average $\frac{1}{n}\sum_{j=1}^n Y_j$ (recall that $n=\lfloor N/k\rfloor$). The term $\|m_n - \mu\|$ is the main contribution to the bias of the estimator $\wh \mu_N$ and is controlled by the size of the block $n$, while $\|\wh\mu_N - m_n\|$ is the stochastic error that depends on the number of blocks $k$. 
We show that under various conditions encoded by the classes $\m P_j, \ j\in\{1,2,3\}$, the ``bias'' admits a dimension-free upper bound of the form $\sqrt{\|\Sigma_Y\|}\sqrt{\frac{k}{N}}$ while 
\ben{
\label{eq:bound2}
\|\wh \mu_N - m_n \|\lesssim \sqrt{\frac{\tr(\Sigma_Y)}{N}} + \sqrt{\Delta\,\tr(\Sigma_Y)} \sqrt{\frac{s}{N}} + \sqrt{\frac{\tr(\Sigma_Y)}{N}}\frac{s}{\sqrt k}
}
with probability at least $1-4e^{-s}$ for all $s\lesssim k$. The combination of \eqref{eq:bound1} and \eqref{eq:bound2} yields the desired inequality.

Our second main result, stated below, is a quadratic growth condition which ensures that any sub-optimality guarantees for the objective function $F(z;y_1,\ldots,y_k)$ translate into the corresponding bounds for the numerical approximation to the geometric median. 

Given a collection of points $y_1,\ldots,y_k\in \mb R^d$ and a positive integer $p$, set $m_p=\frac{1}{k}\sum_{j=1}^k\| y_j\|^p$. 
\begin{theorem}
\label{thm:qgc}
Let $y_1,\ldots,y_k\in \mb R^d$ be such that $\sum_{j=1}^k y_j=0$. Moreover, assume that the matrix 
$\widehat{\Sigma}=\frac{1}{k}\sum_{j=1}^k y_j y_j^T$ satisfies the condition
\[
a:=\sum_{j=2}^d\lambda_j(\widehat{\Sigma})>0
\]
where $\lambda_j(\widehat{\Sigma})$ are the eigenvalues of $\widehat{\Sigma}$ listed in non-increasing order. Then for all $z\in\mb{R}^d$, 
\[
F(z)-F(\wh m)\geq\frac{1}{2} \frac{a\|z-\wh m\|^2}{b^2(\|z-\wh m\|+b)}
\]
where 
\[
b = \frac{20m_1^3+6m_1m_2+m_3}{a}.
\]
\end{theorem}
Note that we do not make any assumptions on the nature of the data $y_1,\ldots,y_k$. 
\begin{remark}
\label{rem:qgcrate}
Observe that whenever $\|z-\wh m\|$ is small, the leading term in the lower bound is $a/2b^3 \|z-\wh m\|^2$. 
If the points are drawn from the uniform distribution on a sphere of radius $\sqrt{d}$, the factor $a/2b^3$ scales like $d^{-1/2}$. Numerical experiments in Section \ref{sec:exp} verify that this rate of dimensional dependence is asymptotically sharp. 
To see that the bound is of correct form in general, consider the collection of points in $R^2$ given by
\[
\begin{pmatrix}
1\\
0
\end{pmatrix}, \begin{pmatrix}
-1\\
0
\end{pmatrix}, \begin{pmatrix}
0\\
1
\end{pmatrix}, \begin{pmatrix}
0\\
-1
\end{pmatrix}
\]
The geometric median of these data is the origin, but the function $F(z)$ restricted to the $x$ axis is
\[
F((x,0)) = 4 + 2 \frac{x^2}{\sqrt{x^2+1}+1},
\]
indicating that local quadratic growth bound is optimal in general.
\end{remark}
\noindent Theorem \ref{thm:qgc} immediately implies the following global error bound.
\begin{corollary}\label{cor:GEB}
Under the assumptions of Theorem \ref{thm:qgc}, for all $z\in\mb{R}^d$, 
\[
\Vert \nabla F(z)\Vert \geq\frac{1}{2} \frac{a\|z-\wh m\|}{b^2(\|z-\wh m\|+b)}.
\]
\end{corollary}
This bound provides a test for early termination given any iterative method. In particular, the right-hand side is monotone increasing, hence
\[
\Vert \nabla F(z)\Vert < \frac{1}{2} \frac{a\varepsilon}{b^2(\varepsilon+b)},
\]
is a sufficient condition for the inequality $\|z-\wh m\|<\varepsilon$ to hold.

\section{Statistical error bounds}
\label{sec:statistical}

In this section, we develop the technical background needed prove Theorem \ref{th:main}. The theorem itself is proved in subsection \ref{sec:mom}.

\subsection{Preliminaries: small ball probabilities}
\label{ssec:sball}

For a centered random vector $Z\in \mb R^d$ with a distribution that is absolutely continuous with respect to the Lebesgue measure, let $M(Z)$ denote the sup-norm of the density $p_Z$ of $Z$. The following ``small-ball'' inequality is immediate: for any $z\in \mb R^d$ and $R>0$,
\be{
\pr{\|Z - z\| \leq R} \leq M(Z) V_d(R)
}
where $V_d(R) = \frac{\l(\sqrt{\pi} R\r)^d}{\Gamma(d/2+1)}$ is the volume of a ball $B(R)$ of radius $R$ in $\mb R^d$. 
Assuming that the covariance matrix $\Sigma_Z = \mb E(Z-\mb EZ)(Z-\mb EZ)^T$ exists (note that it must be non-degenerate for the density $p_Z$ to be well-defined), it is easy to see using the change-of-variables formula that $M(Z) = \frac{M\l(\Sigma_Z^{-1/2}Z\r)}{\sqrt{\det(\Sigma_Z)}}$, hence
\ben{
\label{eq:small-ball0}
\pr{\|Z - z\| \leq R} \leq M\l(\Sigma_Z^{-1/2} Z\r) \frac{V_d(R)}{\sqrt{\det(\Sigma_Z)}}.
}
The advantage of the latter expression is that the quantity $M\l(\Sigma_Z^{-1/2}Z\r)$ is invariant with respect to the affine transformations of $Z$. 
Let us also recall that $V_d(R)$ satisfies the following inequalities for some absolute positive constants $c_1$ and $c_2$: 
\ben{
\label{eq:vol-bound}
\frac{c_1}{\sqrt d} \l( \frac{\sqrt{2\pi e}\,R}{\sqrt{d}}\r)^d \leq V_d(R) \leq \frac{c_2}{\sqrt d} \l( \frac{\sqrt{2\pi e}\,R}{\sqrt{d}}\r)^d.
}
For special classes of distributions, better estimates for the small ball probabilities are available. Next, we will recall several results in this direction.
\begin{theorem}[Theorem 4 in \cite{latala2005small}]  
\label{th:latala}
Let $Z$ have multivariate normal distribution $N(0,\Sigma)$ and let $M$ be the median of $\|Z\|$. Then for all $x\in \mb R^d$,
\[
\pr{\|Z - x\| \leq tM}\leq \frac{1}{2}\l(2t\r)^{\frac{M^2}{4\|\Sigma\|}}.
\]
\end{theorem}
It is helpful to recall that $c_1\sqrt{\tr(\Sigma)}\leq M \leq c_2\sqrt{\tr(\Sigma)}$ for absolute constants $0<c_1<c_2<\infty$, implying that the size of small balls is essentially controlled by the effective rank $r(\Sigma)$. 
A more general result, stated below, is due to Rudelson and Vershynin \cite{rudelson2015small}. 
\begin{theorem}[Theorem 1.5 in \cite{rudelson2015small}]
\label{th:vershynin}
Assume that $Z\in \mb R^d$ is given by a linear transformation $Z = AX$ where $X=(X_1,\ldots,X_k)\in \mb R^{k}$ is a centered random vector with independent coordinates such that the covariance matrix $\Sigma_X = I_{k}$ and $M_0:=\max_{j=1,\ldots,k} M(X_j) < \infty$. 
Then for any $\eps>0$, there exists a positive constant $C_\eps$ such that for all $x\in \mb R^d$ and $t>0$, 
\[
\pr{\|Z - x\|\leq t\sqrt{\tr(\Sigma_Z)}}\leq \l( C_\eps M_0 t\r)^{(1-\eps)\wt r(\Sigma_Z)},
\] 
where $\Sigma_Z = AA^T$ and $\wt r(\Sigma_Z) = \l\lfloor \frac{\tr(\Sigma_Z)}{\|\Sigma_Z\|}\r\rfloor = \lfloor r(\Sigma_Z)\rfloor$. 
\end{theorem}

In the following sections, we will be especially interested in the small ball probabilities associated with $Z_n = \frac{1}{\sqrt{n}} \sum_{j=1}^n (Y_j - \E Y_j)$ where $Y_1,\ldots,Y_n$ are i.i.d. copies of a random vector $Y$ with covariance matrix $\Sigma_Y$. To make the inequality \eqref{eq:small-ball0} useful, we need a non-asymptotic estimate for $M\l(\Sigma_Y^{-1/2}Z_n\r)$. 

To this end, we will rely on two facts. The first is the  generalization of Rogozin's inequality proved by Juvskevivcius and Lee \cite{juvskevivcius2015small}: let $U_1,\ldots,U_n$ be i.i.d. copies of a random vector $U$ with uniform distribution over a ball centered at the origin and with radius $R_U$ such that $M(U) = M\l(\Sigma_Y^{-1/2}Y\r)$. Then 
\ben{
\label{eq:rogozin}
M\l(\Sigma^{-1/2}_YZ_n\r) \leq M\l(\frac{1}{\sqrt n}\sum_{j=1}^n U_j \r).
}
The second estimate, established by Madiman et al. \cite{madiman2017rogozin}, page 17, states that 
\be{
M\l(\frac{1}{\sqrt n}\sum_{j=1}^n \wt U_j\r) \leq c(d):=\frac{(1+d/2)^{d/2}}{\Gamma(1+d/2)}
}
where $\wt U_1,\ldots,\wt U_n$ are i.i.d. with uniform distribution over a ball in $\mb R^d$ of unit volume. 
The definition of $R_U$ yields that $\vol{B\l(R_U\cdot M^{1/d}\l(\Sigma^{-1/2}X\r)\r)} = 1$, 
hence 
\[
M\l( \frac{M^{1/d}\l(\Sigma_Y^{-1/2}Y\r)}{\sqrt n}\sum_{j=1}^n U_j \r) \leq c(d).
\] 
As $M(cY) = c^{-d}M(Y)$ for any random vector $Y\in \mb R^d$, we conclude using \eqref{eq:rogozin} that 
\[
M\l(\Sigma_Y^{-1/2}Z_n\r) \leq M\l(\Sigma_Y^{-1/2}X\r) \frac{(1+d/2)^{d/2}}{\Gamma(1+d/2)}.
\]
Employing the inequality $\Gamma(1+d/2) \geq \sqrt{2\pi d/2}\l( \frac{d}{2e}\r)^{d/2}$, we get a simple bound 
\ben{
\label{eq:density-bound}
M\l(\Sigma_Y^{-1/2}Z_n\r) \leq M\l(\Sigma_Y^{-1/2}Y\r)(2e)^{d/2}
}
and a small ball estimate
\ben{
\label{eq:small-ball}
\pr{\|Z_n - z\| \leq R} \leq c_2 \frac{M\l(\Sigma_Y^{-1/2}Y\r)}{\sqrt{\det(\Sigma_Y)}} \l( \frac{2e\sqrt{\pi}\,R}{\sqrt{d}}\r)^d.
}

\subsection{Upper bounds for the difference between the mean and the median}
\label{section:bias}

In this section, for the ease of notation we will assume that $Y\in \mb R^d$ is centered and that $m$ is the geometric median of $P_Y$. Our goal is to estimate the distance between the mean and the median (which equals $\|m\|$ under our assumptions), hence we will exclude the trivial case $m = 0$. We are especially interested in the situation when the size of $\|m\|$ is independent of or is weakly dependent on the ambient dimension $d$. 
\begin{theorem}
\label{th:bias}
Assume that the distribution of $Y$ is absolutely continuous with respect to Lebesgue measure on some linear subspace of $\mb R^d$. Then
\[
\|m\|\leq \min\l( \sqrt{\tr(\Sigma_Y)}, \sqrt{\|\Sigma_Y\|}  \, \frac{\mb E^{1/2} \|Y-m\|^{-2}}{\mb E \|Y-m\|^{-1}}\r).
\]
\end{theorem}
\begin{proof}
The first part of the bound is straightforward: indeed, since $m$ minimizes the function $z\mapsto \mb E\|Y-z\|$,
\be{
\|m\| = \| m - \mb EY \| \leq \mb E\|Y-m\| \leq \mb E\|Y\| \leq \mb E^{1/2} \|Y\|^2. 
}
To deduce the second inequality, note that under the stated assumptions the median $m$ satisfies the equation $\mb E \frac{Y-m}{\|Y-m\|}=0$, which implies that 
$m = \l(\mb E\frac{1}{\|Y-m\|} \r)^{-1}\mb E\frac{Y}{\|Y-m\|}$. 
Therefore, for any unit vector $u$,
\be{
\dotp{m}{u} =  \l(\mb E\frac{1}{\|Y-m\|} \r)^{-1} \mb E\frac{\dotp{Y}{u}}{\|Y-m\|} \leq \frac{\mb E^{1/2} \|Y-m\|^{-2}}{\mb E \|Y-m\|^{-1}} \mb E^{1/2} \dotp{Y}{u}^2,
}
implying that $\|m\|\leq \sqrt{\|\Sigma_Y\|} \cdot \frac{\mb E^{1/2} \|Y-m\|^{-2}}{\mb E \|Y-m\|^{-1}}$.
\end{proof}
The inequality $\|m\|\leq \sqrt{\tr(\Sigma_Y)}$ is useful when the effective rank $\mathrm{r}\l(\Sigma_Y\r)$ is small. When $\mathrm{r}(\Sigma_Y)$ is large, it is often possible to find a bound for the ratio of negative moments. This problem will be discussed in the following section.

\subsection{Equivalence of the negative moments of the norm}
\label{section:moments}

In view of the inequality stated in Theorem \ref{th:bias}, it is interesting to understand when the ratio $\frac{\mb E^{1/2} \|Y-m\|^{-2}}{\mb E \|Y-m\|^{-1}}$ of negative moments is ``small,'' in particular, when it does not depend on the ambient dimension. We will present several sufficient conditions in this section that cover many typical situations. We state the examples in the order of increasing generality: (a) the case of Gaussian random vectors; (b) the case of linear transformations of a vector with absolutely continuous independent coordinates and (c) the case of absolutely continuous distributions with bounded density. 

\begin{lemma}
\label{lemma:normal}
Assume that $Y$ has normal distribution $N(0,\Sigma_Y)$ such that the effective rank of the covariance matrix $r(\Sigma_Y)> 10$. Then $\frac{\mb E^{1/2} \|Y-m\|^{-2}}{\mb E \|Y-m\|^{-1}} \leq C$ for an absolute constant $C$.
\end{lemma}
\begin{proof}
The claim follows from Theorem \ref{th:latala} (see Corollary 1 in \cite{latala2005small}) once we notice that the median $M(\|Y\|)$ of $\|Y\|$ satisfies $M(\|Y\|)\geq 0.08\sqrt{\tr(\Sigma_Y)}$. Indeed, recall that $Y = \Sigma^{1/2}Z$ where $Z$ has standard normal distribution. Therefore, $\mb E\|Y\| =\mb E\sqrt{Z^T \Sigma_Y Z} = \mb E\sqrt{\sum_{j=1}^d \lambda_j(\Sigma_Y)Z_j^2}=:f(\lambda_1,\ldots,\lambda_d)$. Observe that the function $f$ is concave, hence its minimum in the set 
\[
\l\{ (\lambda_1,\ldots,\lambda_d): \ \lambda_j\geq 0 \ \forall j, \ \sum_{j=1}^d \lambda_j = \tr(\Sigma_Y) \r\}
\]
is achieved at an extreme point $(\tr(\Sigma_Y),0,\ldots,0)$, implying that $\mb E\|Y\| \geq \sqrt{\tr(\Sigma_Y)}\sqrt{\frac2\pi}$. 
It remains to apply Paley-Zygmund inequality to deduce that 
\be{
\pr{\|Y\|\geq t \sqrt{\frac2\pi}\sqrt{\tr(\Sigma_Y)}}  \geq \pr{\|Y\|\geq t \mb E\|Y\|} 
\geq (1-t)^2 \frac{\l(\mb E\|Y\|\r)^2}{\mb E\|Y\|^2} \geq (1-t)^2 \frac{2}{\pi} 
}
which equals $0.5$ for $t = 1 - \sqrt{\pi}/2 >0.11$, and the claim follows. To apply Corollary 1 in \cite{latala2005small}, we require that $\frac{M^2(\|Y\|)}{4\|\Sigma_Y\|}>2$, which holds in view of the previous bound whenever $r(\Sigma_Y)>10$.
\end{proof}

Next, we show that the equivalence of negative moments holds for a larger class of distributions given by linear transformations of a vector with independent coordinates. This class, denoted $\m P_1$, was formally defined in section \ref{sec:main}. 
Since any multivariate normal vector is a linear transformation of the standard normal distribution, Lemma \ref{lemma:independent} below also implies a version of Lemma \ref{lemma:normal}. Recall that $M(Y)$ stands for the sup-norm of the probability density function of a random vector $Y$. 
\begin{lemma}
\label{lemma:independent}
Assume that $Y\in \mb R^d$ has distribution $P_Y$ that belongs to the class $\m P_1(M_0,K)$. Moreover, suppose that the effective rank $r(AA^T)\geq 4$. Then
\[
\frac{\mb E^{1/2} \|Y-m\|^{-2}}{\mb E \|Y-m\|^{-1}}\leq CM_0
\]
for an absolute constant $C>0$.
\end{lemma}
\begin{proof}
Note that $\Sigma_Y = A A^T$. Therefore, 
\[
\l(\mb E \|Y-m\|^{-1}\r)^{-1} \leq \mb E\|Y-m\| \leq \mb E \|Y\| \leq \sqrt{\tr\l(\Sigma_Y\r)}
\] 
in view of Jensen's and Cauchy-Schwarz inequalities. Next, we will prove a general upper bound for $\mb E\|Y-x\|^{-q}$. 
To this end, we will use Theorem 1.5 from the work by Rudelson and Vershynin \cite{rudelson2015small} which states that for any $\eps>0$, there exists a positive constant $C_\eps$ such that for all $x\in \mb R^d$ and $t>0$, $\pr{\|Y-x\|\leq t\sqrt{\tr(\Sigma_Y)}}\leq \l( C_\eps M_0 t\r)^{(1-\eps)\wt r(\Sigma_Y)}$, where $\wt r(\Sigma_Y) = \l\lfloor \frac{\tr(\Sigma_Y)}{\|\Sigma_Y\|}\r\rfloor = \lfloor r(\Sigma_Y)\rfloor$. 
Employing this ``small ball'' bound and letting $r:=\wt r(\Sigma_Y)$ for brevity, we deduce that for any $\delta>0$ and $q<r$,
\ml{
\mb E\|Y - x\|^{-q} = \int_0^\infty \pr{\|Y-x\|\leq t^{1/q}}\frac{dt}{t^2} 
\\
= \frac{q}{\l(\tr(\Sigma_Y)\r)^{k/2}} \l( \int_{1/\delta}^\infty \frac{ds}{s^{q+1}} + \int_0^{1/\delta} (C_\eps M_0 s)^{r(1-\eps)}\frac{ds}{s^{q+1}} \r).
}
Choosing $\delta$ to make the sum above small (e.g. $\delta = C_\eps M_0 \l(\frac{q}{r(1-\eps)-q}\r)^{1/r(1-\eps)}$), it is easy to deduce the inequality
\[
\mb E\|Y - x\|^{-q}\leq \frac{2(C_\eps M_0)^{q}}{\l(\tr(\Sigma_Y)\r)^{q/2}} \l( \frac{q}{r(1-\eps) - q}\r)^{q/r(1-\eps)}.
\]
If $\eps = \frac{r-q-1/2}{r}$, then $\frac{q}{r(1-\eps)-q} = 2q$ and $\frac{q}{r(1-\eps)}\leq 1$. 
For small values of $q$, say, $q\leq r/2$, this choice of $\eps$ entails the inequality $\eps>\frac{r-1}{2r}\geq \frac{3}{8}$ for $r\geq 4$, so that $C_\eps$ can be treated as an absolute constant. 
The claim of the lemma corresponds to the case $q=2$. 
\end{proof}
\noindent Finally, we discuss the most general situation of absolutely continuous distributions. 
\begin{lemma}
\label{lemma:trace}
Assume $Y\in \mb R^d$ has distribution $P_Y$ that belongs to the class $\mathcal{P}_2(k,M_0,K,R)$.  Then for any $x$ in the range $L$ of $\Sigma_Y$ and $q<k=\dim(L)$, 
\[
\mb E\| Y - x\|^{-q}\leq c(q) \frac{M\l(\Sigma_Y^{-1/2}Y\r)^{q/k}}{\l( k\cdot\det^{1/k}(\Sigma_Y)\r)^{q/2}}
\]
for some constant $c(q)>0$.
\end{lemma}
\begin{proof}
The proof is similar to the argument behind Lemma \ref{lemma:independent}. Note that for any $\delta>0$
\begin{align*}
\mb E\| Y - x\|^{-q} 
&=\int_0^\infty \P(\| Y - x\|\leq t^{1/q})t^{-2}\:dt\\
&\leq \int_{1/\delta}^\infty t^{-2}\:dt + \int_0^{1/\delta}\P(\| Y - x\|\leq t^{1/q})t^{-2}\:dt\\
&= \delta + \int_0^{1/\delta}\P(\| Y - x\|\leq t^{1/q})t^{-2}\:dt\\
&\leq \delta + c_2 \frac{M\l(\Sigma^{-1/2}Y\r)}{\sqrt{\det(\Sigma_Y)}}\int_0^{1/\delta}  \l( \frac{\sqrt{2\pi e}\,t^{1/q}}{\sqrt{d}}\r)^k \frac{dt}{t^2}
\end{align*}
in view of \eqref{eq:small-ball}. For the choice of $\delta = c_3(q)\frac{M\l(\Sigma_Y^{-1/2}Y\r)^{q/k}}{\l( k \det^{1/k}(\Sigma_Y)\r)^{q/2}}$, the latter is bounded by $c_4(q)\frac{M\l(\Sigma_Y^{-1/2}Y\r)^{q/k}}{\l( k \cdot\det^{1/k}(\Sigma_Y)\r)^{q/2}}$ for some constant $c_4>0$ that depends only on $q$.
\end{proof}

Since $\l(\mb E \|Y-m\|^{-1}\r)^{-1} \leq \sqrt{\tr\l(\Sigma_Y\r)}$, we immediately get from the previous result that whenever $k\geq 3$, then for some absolute constant $C>0$
\ben{
\label{eq:moments-ratio}
\frac{\mb E^{1/2} \|Y-m\|^{-2}}{\mb E \|Y-m\|^{-1}} \leq CM^{1/k}\l(\Sigma_Y^{-1/2}Y\r) \sqrt{\frac{\tr(\Sigma_Y)}{ k \cdot\det^{1/k}(\Sigma_Y)}}.
}

Lemma \ref{lemma:trace} is robust to small perturbations: for example, assume that $\Sigma_Y = \lambda\sum_{j=1}^k e_j e_j^T + \delta I_d$ where $d\cdot\delta \leq Ck \cdot\lambda$. In this case, $\frac{\tr(\Sigma_Y)}{ d \cdot\det^{1/d}(\Sigma_Y)}$ can be very large, and direct application of Lemma \ref{lemma:trace} yields a suboptimal bound.  
The following simple observation often yields a better result: for any linear subspace $H$ of $\mb R^d$,
\ben{
\label{eq:projection}
\mb E \| Y - x\|^{-q} \leq \mb E \|\Pi_H(Y - x)\|^{-q},
}
where $\Pi_H(\cdot)$ stands for the orthogonal projection onto $H$. We formalize this observation in the following statement.
\begin{lemma}
\label{lemma:projection}
Assume that $Y=X+\xi \in \mb R^d$ has distribution $P_Y$ that belongs to the class $\m P_3(k,M_0,K,R,h)$ and that $k\geq 3$. 
Then for any $x\in \mb R^d$, 
\[
\frac{\mb E^{1/2} \|Y-x\|^{-2}}{\mb E \|Y-x\|^{-1}}\leq C(1+h)M^{1/k}\l( \Sigma_X^{-1/2} X\r) \sqrt{\frac{\tr(\Sigma_X)}{k\det^{1/k}(\Sigma_X)}}.
\]
\end{lemma}
\begin{proof}
Let $H$ be the range of $\Sigma_X$, where, according to the definition of the class $\m P_3$, $Y = X+\xi$. The dimension of $H$ equals $k$ by assumption. Employing the previously stated observation \eqref{eq:projection}, we deduce that 
\al{
\mb E \| Y - x\|^{-2} 
\leq \mb E \|\Pi_H(Y - x)\|^{-2} &= \mb E \| (X + \Pi_H \xi) - \Pi_H x\|^{-2}
\\
&\leq c \frac{M^{2/k}(\wt Y)}{k},
}
where $\wt Y = X + \Pi_H \xi$. It remains to note that $M(\wt Y)\leq M(X)$ by the elementary properties of the convolution operator, and that $M(X) = \frac{M\l( \Sigma_X^{-1/2} X\r)}{\sqrt{\det(\Sigma_X)}}$.
\end{proof}
In the case when $\Sigma_X = \lambda\sum_{j=1}^k e_j e_j^T$, $\Sigma_\xi =  \delta I_d$ and $d\cdot\delta \leq Ck \cdot\lambda$, the previous result yields that the ratio of moments is at most $O(1)M^{1/k}\l( \Sigma_X^{-1/2} X\r)$. 

\begin{remark}
It should be noted that there exist examples where estimates based on the ratios of the arithmetic and geometric means provide only crude bounds: for instance, if $\lambda_j(\Sigma_Y) = \frac{1}{m+j}, \ j=1,\ldots,d$ for some positive integer $m$, then $\frac{\sum_{j=1}^d \lambda_j}{\max_{j\geq 1} j \l(\prod_{i=1}^j \lambda_i\r)^{1/j}}$ can be made arbitrary large by varying $m$ and $d$ (more specifically, it is large when $m/d$ is large). However, under additional assumptions on the distribution (e.g. in the framework of Lemmas \ref{lemma:independent} - \ref{lemma:projection}), better bounds become possible.
\end{remark}
		
\subsection{The geometric median: bounds for the stochastic error}
\label{section:stochastic}

Our goal in this section is to establish high-confidence deviation bounds for the distance between the empirical geometric median and its population counterpart. 
Denote 
\[
D_Y := \mb E\l[\frac{(Y - m)}{\| Y - m\|} \frac{(Y - m)^T}{\| Y - m\|} \r], \quad \Delta = \l\| D_Y \r\|.
\] 
Note that $\tr\l( \mb E\l[\frac{(Y - m)}{\| Y - m\|} \frac{(Y - m)^T}{\| Y - m\|} \r] \r) = 1$. Therefore, if the random vector $Y$ is sufficiently ``spread out,'' we expect that $\Delta$ will be small. To get a rigorous bound supporting this intuition, we will assume that $Y$ satisfies the following conditions:

\begin{enumerate}
\item[(a)] Condition \eqref{eq:moment-equiv} is satisfied. We are especially interested in the situation when $K$ is a constant that does not depend on the ambient dimension $d$.
\item[(b)] $\mb E^{(q-2)/q} \| Y-m \|^{-\frac{q}{2(q-2)}} \leq \frac{C(q)}{\tr(\Sigma_Y)}$.
\end{enumerate}

When (a) and (b) hold, H\"older's inequality implies that
\begin{align}
\Delta = \sup_{\|u\|=1} \mb E \frac{\dotp{Y - m}{u}^2}{\| Y - m\|^2} 
\leq \sup_{\|u\|=1} \mb E^{2/q} \l| \dotp{Y - m}{u}\r|^q \mb E^{(q-2)/q} \| Y-m \|^{-\frac{q}{2(q-2)}} 
\\
\label{eq:sign-covariance}
\leq K\,C(q)\frac{\|\Sigma_Y\|}{\tr(\Sigma_Y)} = \frac{K\,C(q)}{r(\Sigma_Y)}.
\end{align}
Moment equivalence condition (a) has been discussed in detail in section in remark \ref{remark:key}. Condition (b) holds for the classes of distributions discussed in section \ref{section:moments} when the effective rank of 
$\Sigma_Y$ is sufficiently large relative to $\frac{q}{q-2}$. For instance, it holds for linear transformations of random vectors with independent coordinates as well as for random vectors with ``well-conditioned'' covariance matrices, in a sense that the geometric mean of their eigenvalues is equivalent to the arithmetic mean. 
We conclude that for large classes of distributions, $\tr(\Sigma_Y)\Delta \asymp \|\Sigma_Y\|$: indeed,  if $r(\Sigma_Y)$ is small, then it follows since $\Delta \leq 1$, and if $r(\Sigma_Y)$ is large, it follows from the previous discussion. We are ready to state the main result of this section.

\begin{theorem}
\label{th:stochastic}
Let $m:=m(P_Y)$ be the geometric median associated with the distribution $P_Y$ and $\wh m$ - its empirical counterpart based on an i.i.d. sample $Y_1,\ldots,Y_k$ from $P_Y$. Assume that $\Delta < 1$ and that 
\ben{
\label{eq:tech-assump}
\mb E^{1/2}\frac{1}{\|Y - m\|^2} \l(\sqrt{\frac{\tr(\Sigma_Y)}{k}} + \sqrt{\Delta\,\tr(\Sigma_Y)} \sqrt{\frac{s}{k}} + \sqrt{\tr(\Sigma_Y)}\frac{s}{k}\r)<c_1(\Delta).
}
Then 
\[
\|\wh m - m\|\leq K(\Delta)\l(\sqrt{\frac{\tr(\Sigma_Y)}{k}} + \sqrt{\Delta\,\tr(\Sigma_Y)} \sqrt{\frac{s}{k}} + \sqrt{\tr(\Sigma_Y)}\frac{s}{k}\r)
\]
that holds with probability at least $1-2e^{-s}-2e^{-k/4}$ for $s\leq c_2(\Delta)k$.
\end{theorem}
\begin{remark}
\label{remark:stochastic} \.\\

\begin{enumerate}
\item 
In view of the discussion preceding the theorem, we are mainly interested in the situation when $r(\Sigma_Y)$ is bounded from below by a sufficiently large absolute constant and when $\Delta$ is not too close to $1$, e.g. $\Delta\leq 1/2$.
\item Assumption \eqref{eq:tech-assump} is rather mild: indeed, we showed in section \ref{section:moments} that in many common situations, $\mb E^{1/2}\frac{1}{\|Y - m\|^2}\asymp \frac{1}{\sqrt{\tr(\Sigma_Y)}}$.
\item Note that for $s\lesssim \sqrt{k}$, the bound of the theorem yields deviations guarantees of sub-Gaussian type for the median $\wh m$: namely, 
\[
\|\wh m - m\|\leq K(\Delta)\l(\sqrt{\frac{\tr(\Sigma_Y)}{k}} + \sqrt{\Delta\,\tr(\Sigma_Y)} \sqrt{\frac{s}{k}} \r)
\]
with probability at least $1-4e^{-s}$. Despite the fact that asymptotic properties of the geometric median have been well-understood, we believe that this bound is new. 
\end{enumerate}
\end{remark}

\begin{proof}
Recall that, in view of Theorem 3.1 in \cite{minsker2015geometric}, $\|\wh m - m\| \leq 2\sqrt{\tr(\Sigma_Y)}$ on event $\m E$ of probability at least $1-e^{-k/4}$ (it suffices to take $p=1/8$ and $\alpha = 5/12$ in the aforementioned result). In what follows, we will assume that event $\m E$ occurs. 
Define $\wh u := \frac{m - \wh m}{\|m-\wh m\|}$ (for absolutely continuous distributions, $\wh m \ne m$ with probability 1, so $\wh u$ is well-defined) and 
\[
G_k(s) := \frac{1}{k}\sum_{j=1}^k \| m + s\wh u - Y_j\|.
\]
Then $G_k(s)$ is convex, achieves its minimum at $\wh s = \|\wh m - m\|$, and its derivative $G_k'(s)$ is non-decreasing and satisfies $G_k'(s)\leq 0$ for $s\in [0,\wh s]$. It implies that $\|\wh m - m\| \geq t$ is true only if $G_k'(t)\leq 0$. In view of convexity of $G_k$, 
\[
0\geq G_k'(t) \geq G_k'(0) + \inf_{0\leq z\leq t} G_k''(z)\cdot t,
\]
where $G_k''(z) = \frac{1}{k}\sum_{j=1}^k \frac{1}{\| m + z\wh u - Y_j\|}\l( 1 - \dotp{\frac{m + z\wh u - Y_j}{\| m + z\wh u - Y_j\|}}{\wh u}^2\r)$. Therefore, a necessary condition for the inequality $\|\wh m - m\| \geq t$ to hold is
\[
\frac{1}{k}\sum_{j=1}^k \dotp{\frac{m - Y_j }{\| m-Y_j \|}}{\wh u} \geq 
t \,\inf_{0\leq z\leq t} \frac{1}{k}\sum_{j=1}^k \frac{1}{\| m + z\wh u - Y_j\|}\l( 1 - \dotp{\frac{m + z\wh u - Y_j}{\| m + z\wh u - Y_j\|}}{\wh u}^2\r),
\]
which is possible only if 
\[
\l\|\frac{1}{k}\sum_{j=1}^k \frac{ m - Y_j}{\| m - Y_j \|} \r\| \geq 
t \,\inf_{0\leq z\leq t} \frac{1}{k}\sum_{j=1}^k \frac{1}{\| m + z\wh u - Y_j\|}\l( 1 - \dotp{\frac{m + z\wh u - Y_j}{\| m + z\wh u - Y_j\|}}{\wh u}^2\r).
\]
Next, we will find high confidence bounds for both sides of the inequality above. Note that we can assume that $t\leq 2\sqrt{\tr(\Sigma_Y)}$ on event $\m E$.
\begin{lemma}
\label{lemma:signs}
With probability at least $1-e^{-s}$, 
\ben{
\label{eq:event1}
\l\| \frac{1}{k}\sum_{j=1}^k \frac{Y_j - m}{\| Y_j - m\|} \r\| \leq \frac{2}{\sqrt k} +\sqrt{\Delta}\sqrt{\frac{2s}{k}} + \frac{4s}{3k}.
}
Moreover, if that random vector $\frac{Y-m}{\|Y-m\|}$ has sub-Gaussian distribution, then 
\ben{
\label{eq:event2}
\l\| \frac{1}{k}\sum_{j=1}^k \frac{Y_j - m}{\| Y_j - m\|} \r\| \leq C\l(\frac{1}{\sqrt k} +\sqrt{\Delta}\sqrt{\frac{s}{k}} \r)
}
for an absolute constant $C>0$ and with probability at least $1-e^{-s}$.
\end{lemma}
\begin{proof}
Let $X_k(u) =  \frac{1}{k}\sum_{j=1}^k \dotp{\frac{Y_j - m}{\| Y_j - m\|}}{u}$ and note that $\mb EX_k(u) = 0$ for all $u$. Next, write the norm as 
\[
\sup_{\|u\|=1} X_k(u) = \sup_{\|u\|=1} \frac{1}{k} \sum_{j=1}^k \dotp{\frac{Y_j - m}{\| Y_j - m\|}}{u}.
\]
Bousquet's version of Talagrand's concentration inequality (see \cite{boucheron2013concentration}) yields that 
\[
\sup_{\|u\|=1} X_k(u) \leq 2\mb E \sup_{\|u\|=1} X_k(u) + \sup_{\|u\|=1}\var^{1/2}\l( X_k(u)\r)\sqrt{2s} + \frac{4s}{3k} 
\]
with probability at least $1-e^{-s}$. It remains to notice that 
\[
\mb E \sup_{\|u\|=1} X_k(u) \leq \mb E^{1/2}\l\| \frac{1}{k}\sum_{j=1}^k \frac{Y_j - m}{\| Y_j - m\|} \r\|^2 = \frac{1}{\sqrt k} \mb E^{1/2} \l\|\frac{Y_1 - m}{\| Y_1 - m\|} \r\|^2 = \frac{1}{\sqrt k}
\]
and that $\sup_{\|u\|=1}\var^{1/2}\l( X_k(u)\r) = \frac{1}{\sqrt k}\l\| \mb E \frac{Y_1 - m}{\| Y_1 - m\|} \frac{(Y_1 - m)^T}{\| Y_1 - m\|}\r\|^{1/2} = \sqrt{\frac{\Delta}{k}}$. 
Part (b) of the lemma follows from the standard concentration bound for sub-Gaussian processes (see \cite{dirksen2015tail}) in place of Bousquet's inequality. 
\end{proof}

\begin{lemma}
\label{lemma:hessian-lower}
Let $\tau>$ be a positive constant, and define 
\ml{
\delta:=\delta(k,t,\tau,\Delta;s):=(1+\tau)\l(\sqrt{\Delta}\l( 1+\sqrt{\frac{2s}{k}}\r) + \frac{4}{\sqrt k}\r) 
\\
+ 2(4+1/\tau)t^2 \mb E\frac{1}{\|Y - m\|^2}  + \l( 8 + \frac{4\tau}{3}+\frac{5}{3\tau}\r)\frac{s}{k}
}
If $\delta<1$, then the following inequality holds with probability at least $1-2e^{-k/4}-2e^{-s}$:
\[
\inf_{\|u\|=1, \|m-x\|\leq t} \frac{1}{k}\sum_{j=1}^k \frac{1}{\| Y_j - x\|}\l( 1 - \dotp{\frac{Y_j - x}{\| Y_j - x\|}}{u}^2\r) \geq 
\frac{C\l(\delta \r)}{\sqrt{\tr(\Sigma_Y)}}.
\]
\end{lemma}
\begin{remark}
Since $\Delta<1$ (recall that we are mostly interested in the situation $\Delta\leq 1/2$), there exist $\tau=\tau(\Delta)>0$ and $\eps=\eps(\Delta) > 0$ such that $\delta<1$ whenever 
\[
t < \eps \l( \mb E^{1/2} \|Y_1 - m \|^{-2} \r)^{-1}
\]
and $k$ is sufficiently large; let us again recall that in many typical situations, 
\[
\l(\mb E^{1/2} \|Y_1 - m \|^{-2} \r)^{-1} \asymp \sqrt{\tr(\Sigma_Y)}.
\]
\end{remark}

\begin{proof}
Note that on event $\m E$ that was defined at the start of the proof of the theorem, $\|Y_j - x\|\leq \|Y_j - m \| + 2\sqrt{\tr(\Sigma_Y)}$ for all $j$, hence one easily gets that for any $\kappa>0$,
\ml{
\pr{ \exists J\subset [k]: \ |J|\geq \kappa k \text{ and } \| Y_j - x \|\geq (c(\kappa)+2)\sqrt{\tr(\Sigma_Y)}, \ j\in J} 
\\
\leq \binom{k}{\lfloor \kappa k\rfloor}\l( 2/c(\kappa)^2 \r)^{\lfloor \kappa k\rfloor} 
\leq e^{-k}
}
where $c(\kappa) < (e^{1/\kappa}/\kappa)^{1/2}$. Consequently, on event $\m E_1$ of probability at least $1-e^{-k}$, 
\be{
\label{eq:kappa}
\|Y_j - x\|\leq (c(\kappa)+2)\sqrt{\tr(\Sigma_Y)} \text{ for all } j\in J \text{ such that } |J|\geq (1-\kappa)k.
}
Next, we will find an upper bound for $\frac{1}{k}\sum_{j=1}^k \dotp{\frac{Y_j - x}{\| Y_j - x\|}}{u}^2$ that holds uniformly over $u$. Recall the following elementary inequality that is valid for all vectors $y_1,y_2\in \mb R^d$: 
$\l\| \frac{y_1}{\|y_1\|} - \frac{y_2}{\|y_2\|}\r\| \leq 2\frac{\|y_1-y_2\|}{\max\l(\|y_1\|,\|y_2\|\r)}$. It implies that for all $j, \ 1\leq j\leq k$,
\[
\l| \dotp{\frac{Y_j - x}{\| Y_j - x\|}}{u} - \dotp{\frac{Y_j - m}{\| Y_j - m\|}}{u}\r|^2 \leq 4\min\l(1,\frac{\|x-m\|^2}{\|Y_j-m\|^2}\r)
\]
so that for any $\tau>0$,
\mln{
\label{eq:error-decomposition}
\sup_{\|u\|=1,\|x-m\|\leq t} \frac{1}{k}\sum_{j=1}^k \dotp{\frac{Y_j - x}{\| Y_j - x\|}}{u}^2 
\\
\leq \frac{1+\tau}{k}\sum_{j=1}^k \dotp{\frac{Y_j - m}{\| Y_j - m\|}}{u}^2 + \frac{4+1/\tau}{k}\sum_{j=1}^k \min\l(1,\frac{t^2}{\|Y_j-m\|^2}\r). 
}
The first term in the sum above can be estimated as follows: note that 
\[
\frac{1}{k}\sum_{j=1}^k \dotp{\frac{Y_j - m}{\| Y_j - m\|}}{u}^2 \leq \frac{1}{k}\sum_{j=1}^k \l|\dotp{\frac{Y_j - m}{\| Y_j - m\|}}{u}\r|
\]
and define
\[
Z_k(u) = \frac{1}{k}\sum_{j=1}^k \l|\dotp{\frac{Y_j - m}{\| Y_j - m\|}}{u}\r| - \mb E \l|\dotp{\frac{Y_1 - m}{\| Y_1 - m\|}}{u}\r|.
\] 
Bousquet's version of Talagrand's concentration inequality yields that 
\[
\sup_{\|u\|=1} Z_k(u) \leq 2\mb E \sup_{\|u\|=1} Z_k(u) + \sup_{\|u\|=1}\var^{1/2}\l( Z_k(u)\r)\sqrt{2s} + \frac{4s}{3k} 
\]
with probability at least $1-e^{-s}$. It remains to note that 
$\mb E \l|\dotp{\frac{Y_1 - m}{\| Y_1 - m\|}}{u}\r| \leq \sqrt{\Delta}$ in view of Cauchy-Schwarz inequality, and that 
\ml{
\mb E \sup_{\|u\|=1} Z_k(u) \leq 2\mb E\sup_{\|u\|=1} \frac{1}{k}\sum_{j=1}^k \eps_j \l|\dotp{\frac{Y_j - m}{\| Y_j - m\|}}{u}\r| 
\\
\leq 4 \mb E\sup_{\|u\|=1} \frac{1}{k}\sum_{j=1}^k \dotp{\frac{Y_j - m}{\| Y_j - m\|}}{u} 
= 4\mb E\l\|\frac{1}{k}\sum_{j=1}^k \frac{Y_j - m}{\| Y_j - m\|} \r\| \leq \frac{4}{\sqrt k}
}
in view of the symmetrization and Talagrand's contraction inequalities (see \cite{gine2015mathematical}). To summarize, we showed that with probability at least $1-e^{-s}$, for all unit vectors $u$,
\ben{
\label{eq:first-term}
\frac{1+\tau}{k}\sum_{j=1}^k \l|\dotp{\frac{Y_j - m}{\| Y_j - m\|}}{u}\r| 
\leq (1+\tau)\l(\sqrt{\Delta}\l( 1+\sqrt{\frac{2s}{k}}\r) + \frac{4}{\sqrt k} + \frac{4s}{3k}\r).
}
In view of Bernstein's inequality, the second term in \eqref{eq:error-decomposition} is at most 
\mln{
\label{eq:second-term}
(4+1/\tau)\l(\mb E\min\l(1,\frac{z^2}{\|Y_j-m\|^2}\r) + 2\sqrt{\var\l( \min\l(1,\frac{z^2}{\|Y - m\|^2}\r)\r)}\sqrt{\frac{s}{k}} + \frac{2s}{3k} \r) \\
\leq (4+1/\tau)\l(2 t^2\mb E\frac{1}{\|Y - m\|^2} +  \frac{5s}{3k} \r)
}
with probability at least $1-e^{-s}$. 
Combining \eqref{eq:error-decomposition}, \eqref{eq:first-term}, \eqref{eq:second-term}, we deduce the inequality 
\ml{
\sup_{\|u\|=1,\|x-m\|\leq t} \frac{1}{k}\sum_{j=1}^k \dotp{\frac{Y_j - x}{\| Y_j - x\|}}{u}^2  
\\
\leq \delta(k,t,\tau,\Delta;s)
:= (1+\tau)\l(\sqrt{\Delta}\l( 1+\sqrt{\frac{2s}{k}}\r) + \frac{4}{\sqrt k}\r) 
\\
+ 2(4+1/\tau)t^2\mb E\frac{1}{\|Y  - m\|^2}  + \l( 8 + \frac{4\tau}{3}+\frac{5}{3\tau}\r)\frac{s}{k}
}
that holds with probability at least $1-2e^{-s}$. If $\delta(k,t,\tau,\Delta;s) <1$, then
\[
\l| \l\{j: \ \dotp{\frac{Y_j - x}{\| Y_j - x\|}}{u}^2  \geq \delta^{1/2}(k,t,\tau,\Delta;s) \r\}\r| \leq \delta^{1/2}(k,t,\tau,\Delta;s) k
\] 
uniformly over all $\|u\|=1$ and $\|x - m\|\leq t$ with probability at least $1-2e^{-s}$. 
Now we set $\kappa:=\frac{1-\delta^{1/2}(k,t,\tau,\Delta;s)}{2}$ in \eqref{eq:kappa} and deduce that for all $u$, there exists a subset $J$ of cardinality at least $\kappa k$ such that $\frac{1}{\|Y_j - x\|}\geq \frac{1}{C(\kappa)\sqrt{\tr(\Sigma_Y)}}$ and 
$\dotp{\frac{Y_j - x}{\| Y_j - x\|}}{u}^2 < \delta^{1/2}(k,t,\tau,\Delta;s) < 1$ for all $j\in J$. Consequently, 
\be{
\inf_{\|u\|=1, \|m-x\|\leq z} \frac{1}{k}\sum_{j=1}^k \frac{1}{\| Y_j - x\|}\l( 1 - \dotp{\frac{Y_j - x}{\| Y_j - x\|}}{u}^2\r) \geq 
\frac{C\l(\kappa \r)}{\sqrt{\tr(\Sigma_Y)}}
}
with probability at least $1-2e^{-k/4}-2e^{-s}$, where $C\l(\kappa\r) \to \infty$ as $\kappa\to 0$.  
\end{proof}
\noindent To complete the proof of the theorem, choose 
\[
t = \wh t:=K\l(\sqrt{\frac{\tr(\Sigma_Y)}{k}} + \sqrt{\Delta\,\tr(\Sigma_Y)} \sqrt{\frac{s}{k}} + \sqrt{\tr(\Sigma_Y)} \frac{s}{k}\r)
\] 
where the constant $K$ is sufficiently large (the specific requirement for the size of $K$ is given below). If $k\geq k_0(\Delta)$ is  large enough, $s\leq c_1(\Delta)k$ and $\wh t \mb E^{1/2}\|Y-m\|^{-2} \leq c_2(\Delta)$, then $\delta(k,\wh t,\tau,\Delta)<1$, implying that the results of Lemmas \ref{lemma:signs} and \ref{lemma:hessian-lower} hold with $t=\wh t$ on event $\m E_2$ of probability at least $1-2e^{-k/4}-2e^{-s}$. If $\|\wh m\|\geq \wh t$, then the following inequality must hold on $\m E_2$:
\[
\wh t \leq \frac{1}{C'(\Delta)}\l( \sqrt{\frac{\tr(\Sigma_Y)}{k}} + \sqrt{\Delta\,\tr(\Sigma_Y)} \sqrt{\frac{s}{k}}\r).
\]
If $K$ is set so that $K>\frac{1}{C'(\Delta)}$, this yields a contradiction. Finally, the bound for the case when $\frac{Y-m}{\|Y-m\|}$ has sub-Gaussian distribution follows with \eqref{eq:event2} in place of \eqref{eq:event1}.
\end{proof}
	
\subsection{Implications for the median of means estimator}
\label{sec:mom}

In this section we prove Theorem \ref{th:main}. To this end, we will apply Theorems \ref{th:bias} and \ref{th:stochastic} to the distribution $P^{(n)}$ of the average $\frac{1}{n}\sum_{j=1}^n Y_j$ and the sample $\bar Y_1,\ldots,\bar Y_k$, noting that the corresponding covariance matrix satisfies $\Sigma_{\bar Y_1} = \frac{\Sigma}{n}\preceq 2\Sigma\frac{k}{N}$ whenever $k\leq N/2$. In what follows, let $m_n$ denote the geometric median of $P^{(n)}$.

Consider two scenarios: if $r(\Sigma_Y)\leq c\frac{q}{q-2}$, then the inequality \eqref{eq:mom-basic} readily yields the result. 
On the other hand, if $r(\Sigma_Y) > c\frac{q}{q-2}$, then $\Delta\tr(\Sigma_Y)\leq C(q)\|\Sigma_Y\|$ and $\mb E^{1/2}\frac{1}{\|Y-m\|^2}\leq \frac{C'}{\sqrt{\tr(\Sigma_Y)}}$ for a constant $C'$ that depends on the parameters of the class $\m P_j, \ j\in \{1,2,3\}$. It remains to show that the relevant parameters of the distribution $P^{(n)}$ can be controlled by the corresponding parameters of the distribution $P_Y$. 
First, recall the inequality \eqref{eq:density-bound} which implies that 
\[
M^{1/k}\l(\Sigma_{\bar Y_1}^{-1/2}\sqrt{n}\bar Y_1\r) \leq \sqrt{2e} M^{1/k}\l(\Sigma_Y^{-1/2}Y \r).
\] 
Therefore, the ratio 
$\frac{\mb E^{1/2} \|\sqrt{n}(\bar Y_1-m_n)\|^{-2}}{\mb E \| \sqrt{n}(\bar Y_1-m_n)\|^{-1}}$ can be estimated via Lemma \ref{lemma:independent}, Lemma \ref{lemma:trace} or Lemma \ref{lemma:projection} in terms of parameters of the distribution $P_Y$ whenever it belongs to one of the classes $\m P_j,\ j\in\{1,2,3\}$. 
Next, consider the norm of the spatial sign covariance matrix 
\[
\Delta^{(n)} := \l\|\mb E\l[\frac{(\bar Y_1 - m_n)}{\| \bar Y_1 - m_n\|} \frac{(\bar Y_1 - m_n)^T}{\| \bar Y_1 - m_n\|} \r] \r\|.
\] 
In view of the well-known moment bounds (e.g., the Marcinkiewicz-Zygmund type inequality by Rio \cite{rio2009moment}), for any unit vector $u$ and $q>2$, 
\[
\mb E \l|\dotp{\frac{1}{\sqrt n}\sum_{j=1}^n Y_j}{u} \r|^q \leq (q-1)^{q/2} \mb E \l| \dotp{Y_1}{u} \r|^q,
\]
thus the reasoning similar to \eqref{eq:sign-covariance} implies that 
\ben{
\label{eq:sign-cov}
\Delta^{(n)}\leq \frac{KC_1(q)}{r(\Sigma_Y)}
}
whenever $P_Y\in \m P_j, \ j\in\{1,2,3\}$. Therefore, conditions of Theorem \ref{th:stochastic} hold for $k$ large enough, and we deduce that
\mln{
\label{eq:general-bound}
\|\wh \mu_N - \mu\|\leq  \l\| m_n - \mu \r\| + \l\|\wh\mu_N - m_n\r\|
\\
\leq
2 \sqrt{\frac{\|\Sigma_Y\|k}{N}}  \, \frac{\mb E^{1/2} \|\sqrt{n}(\bar Y_1-m_n)\|^{-2}}{\mb E \| \sqrt{n}(\bar Y_1-m_n)\|^{-1}}
+ K(\Delta^{(n)})\l(\sqrt{\frac{\tr(\Sigma_Y)}{N}} + \sqrt{\Delta^{(n)}\,\tr(\Sigma_Y)}\sqrt{\frac{\sqrt k}{N}}    \r)
} 
with probability at least $1-4e^{-\sqrt k}$. The final form of the bound follows once we apply the inequality \eqref{eq:sign-cov} and estimate the ratio of moments via one of the lemmas in section \ref{section:moments}. For instance, if $P_Y\in \m P_1(M_0,K)$, then Lemma \ref{lemma:independent} combined with \eqref{eq:general-bound} implies that 
\ml{
\|\wh \mu_N - \mu\|\leq C\l( M_0  \sqrt{\frac{\|\Sigma_Y\|k}{N}} + \sqrt{\frac{\tr(\Sigma_Y)}{N}} + \sqrt{KC_1(q)}\sqrt{\|\Sigma_Y\|\frac{\sqrt k}{N}} \r) 
\\
\leq C\l( M_0  \sqrt{\frac{\|\Sigma_Y\| k}{N}} + \sqrt{\frac{\tr(\Sigma_Y)}{N}} \r)
}
with probability at least $1-4e^{-\sqrt k}$ whenever $k\geq k_0(M_0,K,q)$. Bounds for the classes $\m P_2$ and $\m P_3$ follow similarly.

\section{Numerical error bounds}
\label{sec:numerical} 


In this section we discuss numerical procedures with provable error bounds based on Theorem \ref{thm:qgc}. 
Theorem \ref{thm:qgc} itself is proved in section \ref{proof:lowlip}.

\noindent For $i=1,\ldots,k$, let $f_i(z) = \Vert z - Y_i\Vert$, and observe that a weak gradient of $f_i(z)$ is given by
\[
\nabla f_i(z) = \left\{\begin{array}{cl}
				0 & z=Y_i\\
				\frac{1}{\Vert z - Y_i\Vert}(z-Y_i) & z\not=Y_i
				\end{array}\right..
\] 
Hence, $\nabla F(z) = \frac{1}{k}\sum_{i=1}^k \nabla f_i(z)$ is a weak gradient of $F$. 

\subsection{Algorithms}
\label{ssec:alg}

Given $\{Y_i\}_{i=1}^k\subset\mb R^d$ and $\delta>0$, we set
\begin{align} 
F_\delta(z) = F_\delta\l(z;\{Y_i\}_{i=1}^k\r) = \frac{1}{k} \sum_{i=1}^k \sqrt{\| z-Y_i\|^2 +\delta^2}\label{def:relax}
\end{align}
We will call this function and the associated minimization program the \emph{Charbonnier relaxation} of the mean norm deviation. The Charbonnier relaxation is smooth, and hence we may perform accelerated gradient descent to approximate its solution. 
\begin{algorithm}
\caption{Accelerated Gradient Descent of the Charbonnier Relaxation}
\label{alg:cragd}
\begin{algorithmic}
\STATE $\varepsilon>0$, $t=0$, $\alpha_0=3/4$, $\{Y_i\}_{i=1}^k\in\mb R^d$, $x^{(0)}=v^{(0)}=\bar Y_k:=k^{-1} \sum_{i=1}^k Y_i$
\WHILE{$\varepsilon/2 \leq \frac{16}{9(t+1)^2}\left[ F_{\varepsilon/2}(\bar Y_k) +\frac{9}{4\varepsilon}F_{\varepsilon/2}(\bar Y_k)^2\right]$}
    \STATE $x^{(t+1)} \gets v^{(t)} - \frac{\varepsilon}{2} \nabla F_{\varepsilon/2}(v^{(t)})$
    \STATE $\alpha_{t+1} \gets \frac{\alpha_t}{2}\left( \sqrt{\alpha_t^2 + 4} - \alpha_t\right)$
    \STATE $v^{(t+1)} \gets x^{(t+1)} + \frac{\alpha_t(1-\alpha_t)}{\alpha_t^2 +\alpha_{t+1}}(x_{t+1}-x_t)$
    \STATE $t\gets t+1$
\ENDWHILE\\
\RETURN $x^{(t)}$
\end{algorithmic}
\end{algorithm}

Using standard results for accelerated gradient descent and simple sub-optimality bounds, we will prove the following estimate in section \ref{ssec:algerr}.

\begin{theorem}
\label{thm:subopt0}
The output of Algorithm \ref{alg:cragd} satisfies $F(x^{(t)}) - F(\wh m) < \varepsilon$.
\end{theorem}

Recall that, in view of Theorem \ref{thm:qgc}, whenever $\|z-\wh m\|$ is small, it behaves like 
\[
\sqrt{F(z)-F(\wh m)}.
\] 
It turn, it implies that the Restarted Gradient Descent algorithm \cite{yang2015rsg} achieves an iteration complexity of order $\mathcal{O}(\varepsilon^{-1}\log(\varepsilon^{-1}))$ for computation of the geometric median. However, our Algorithm \ref{alg:cragd} admits the better iteration complexity $\mathcal{O}(\varepsilon^{-1})$ due to the strong convexity of the Charbonnier relaxation.
	
In practice, we found that the following relaxed Newton's method (Algorithm \ref{alg:crnewt}) provides the best results, especially when combined with the stopping rule described in Corollary \ref{cor:GEB}. 
\begin{algorithm}
\caption{Newton's Method for Successive Charbonnier Relaxation}
\label{alg:crnewt}
\begin{algorithmic}
\REQUIRE $\varepsilon>0$, $t=0$, $\tau=1$, $M>1$, $\{Y_i\}_{i=1}^k\in\mb R^d$, $x^{(0)}=\bar Y_k$
\WHILE{$\frac{a}{2b^2}\leq \Vert F(x^{(t)})\Vert $}
    \STATE $\tau \gets \tau/M$
    \STATE $x^{(t+1)} \gets x^{(t)} - \nabla^2 F_{\tau}(x^{(t)})^{-1}\nabla F_{\tau}(x^{(t)})$
    \STATE $t\gets t+1$
\ENDWHILE
\WHILE{$\varepsilon \leq \frac{2b^3\Vert \nabla F(x^{(t)})\Vert}{a-2b^2 \Vert \nabla F(x^{(t)})\Vert}$}
    \STATE $\tau \gets \tau/M$
    \STATE $x^{(t+1)} \gets x^{(t)} - \nabla^2 F_{\tau}(x^{(t)})^{-1}\nabla F_{\tau}(x^{(t)})$
    \STATE $t\gets t+1$
\ENDWHILE\\
\RETURN $x^{(t)}$
\end{algorithmic}
\end{algorithm}
However, we note that, unlike Algorithm \ref{alg:cragd}, no rigorous analysis for Algorithm \ref{alg:crnewt} is currently available.

\subsection{Proof of Theorem \ref{thm:subopt0}}
\label{ssec:algerr}

To prove Theorem \ref{thm:subopt0}, we first exhibit some simple sub-optimality bounds for the solution to the Charbonnier relaxation. 
Since $F_\delta$ is smooth, the Lipschitz constant for $\nabla F_\delta$ is $\max \|\nabla^2 F_\delta(x)\|$ over $x\in\R^d$. We have the straightforward estimate
\[
\| \nabla^2 F_\delta(x)\| = \left\| \frac{1}{k} \sum_{i=1}^k \frac{1}{\sqrt{\| x-Y_i\|^2 + \delta^2}}\left(I -  \frac{1}{\|x-Y_i\|^2 + \delta^2}(Y_i-x)(Y_i-x)^T\right)\right\| \leq \frac{1}{\delta}
\]
since
\[
\left\| I -  \frac{1}{\|x-Y_i\|^2 + \delta^2}(Y_i-x)(Y_i-x)^T\right\| \leq 1
\]
for all $i=1,\ldots, k$. Thus, we may minimize $F_\delta$ using the ``constant step scheme II'' on page 93 of the book by Nesterov \cite{nesterov2018lectures}. In the notation of \cite{nesterov2018lectures}, we have that $q_f=0$, so we get the following sub-optimality bound.

\begin{lemma}
\label{lem:subopt}
Suppose $\{Y_i\}_{i=1}^k\subset \mathbb{R}^d$, $\delta>0$, and $x_{\ast,\delta}$ minimizes $F_\delta$. Then
\[
F_\delta(x^{(t)}) - F_\delta(x_{\ast,\delta}) \leq \frac{16}{9(t+1)^2}\left[ F_\delta(\bar Y_k) +\frac{9}{8\delta}F_\delta(\bar Y_k)^2\right]
\]
for all $t\geq 0$ given the above algorithm.
\end{lemma}

\begin{proof}\label{proof:subopt}
Theorem 2.2.3 in the book by Nesterov \cite{nesterov2018lectures} holds with $q_f=0$ since 
\[
\alpha_0  = \frac{3}{4} \leq \frac{6}{3+\sqrt{21}}
\]
in accordance with condition (2.2.21) of the book, and gives the bound 
\[
F_\delta(x^{(t)}) - F_\delta(x_{\ast,\delta}) \leq  \frac{4L}{\gamma_0(k+1)^2}\left[ F_\delta(x^{(0)})-F_\delta(x_{\ast,\delta}) + \frac{\gamma_0}{2}\Vert x^{(0)}-x_{\ast,\delta}\Vert^2\right]
\]
with $L=\frac{1}{\delta}$ and $\gamma_0 = \frac{\alpha_0^2}{1-\alpha_0} = \frac{9}{4} L$. The result follows from this bound when we also invoke the inequalities $F_\delta(\bar Y_k)-F_\delta(x_{\ast,\delta}) \leq F_\delta(\bar Y_k)$ and $\| \bar Y_k - x_{\ast,\delta}\| \leq F_\delta(\bar Y)$.
\end{proof}
\noindent We are now ready to prove Theorem \ref{thm:subopt0}.

\begin{proof}[Proof of Theorem \ref{thm:subopt0}]
\label{proof:subopt0}
Using sub-additivity of the square root, we observe that
\[
F(z)\leq F_\delta(z) \leq F(z) + \delta
\]
for all $z\in\mathbb{R}^d$. Consequently, if $\wh m$ minimizes $F$ and $x_{\ast,\delta}$ minimizes $F_\delta$, then
\[
F(\wh m) \leq F(x_{\ast,\delta}) \leq F_\delta(x_{\ast,\delta}) \leq F(\wh m) + \delta.
\]
Therefore, a $\frac{\varepsilon}{2}$ numerical approximation $\tilde x$ to $x_{\ast,\varepsilon/2}$ satisfies
\begin{align*}
F(\tilde{x}) - F(\wh m) & \leq F_\delta(\tilde{x}) - F(\wh m) \\
& = F_\delta(\tilde{x}) - F_\delta(x_{\ast,\varepsilon/2}) + F_\delta(x_{\ast,\varepsilon/2}) - F(\wh m)\\
& \leq \frac{\varepsilon}{2} + \frac{\varepsilon}{2} = \varepsilon.
\end{align*}
With $\tilde{x}=x^{(t)}$ being the output of Algorithm \ref{alg:cragd}, the termination condition and the bound of Lemma \ref{lem:subopt} together imply the result.
\end{proof}

\begin{remark}
Now we have that $F_\delta(x^{(t)}) - F_\delta(x_{\ast,\delta})$ whenever $\frac{4}{3\sqrt{\delta}} \left[f + \frac{9}{8\delta}f^2\right]^{1/2} - 1\leq t$ if $f=F_\delta(\bar Y_k)$. When $\delta=\varepsilon/2$, we see that $t$ is on the order of $1/\varepsilon$. While computer science literature (e.g. \cite{cohen2016geometric}) exhibits better bounds than these, the numerical constants in those algorithms can be impractical. 

\end{remark}

\subsection{Proof of Theorem \ref{thm:qgc}}
\label{proof:lowlip}

Fix $z\in\mb{R}^d$ with $z\not=\wh m$, let $r=\Vert z-\wh m\Vert$, and set $u=\frac{1}{r}(z-\wh m)$. In view of the second fundamental theorem of calculus, 
\begin{align*}
F(z) - F(\wh m) &= \int_0^r \nabla F(\wh m + t u)^T u dt\\
&=\frac{1}{k}\int_0^r \sum_{i=1}^k \frac{1}{\|\wh m-Y_i+ t u\|}(\wh m-Y_i+tu)^Tu dt\\
&=\frac{1}{k}\int_0^r \sum_{i=1}^k \frac{(\wh m-Y_i)^Tu + t}{\sqrt{\|\wh m-Y_i\|^2+2t(\wh m-Y_i)^Tu + t^2}}dt\\
&=\frac{1}{k}\int_0^r \sum_{i=1}^k \frac{\gamma_i c_i + t}{\sqrt{(\gamma_ic_i + t)^2+\gamma_i^2(1-c_i^2)}}dt.
\end{align*}
In this last line, we have set $\gamma_i = \|\wh m - Y_i\|$ and $c_i= \frac{1}{\gamma_i}(\wh m-Y_i)^T u$. By the Cauchy-Schwarz inequality, we have that $c_i^2\leq 1$. If $c_i^2=1$, then 
\[
\frac{\gamma_i c_i + t}{\sqrt{(\gamma_i c_i +t)^2 + \gamma_i^2(1-c_i^2)}}= \text{sgn}(\gamma_i c_i+t)\geq c_i
\]
for all $t> 0$. If $c_i^2<1$, then 
\[
\frac{\gamma_i c_i + t}{\sqrt{(\gamma_i c_i +t)^2 + \gamma_i^2(1-c_i^2)}}=c_i+\int_0^t \frac{\gamma_i^2(1-c_i^2)}{\left[(\gamma_ic_i+s)^2+\gamma_i^2(1-c_i^2)\right]^{3/2}}ds.
\]
Note that $\sum_{i=1}^k c_i = \nabla F(\wh m)^T u=0$ since $\wh m$ is the minimizer. Consequently, we have that
\begin{align*}
F(z) - F(\wh m) &\geq \frac{1}{k}\int_0^r \left(\sum_{i=1}^k c_i + \sum_{i:c_i^2<1} \int_0^t \frac{\gamma_i^2(1-c_i^2)}{\left[(\gamma_ic_i+s)^2+\gamma_i^2(1-c_i^2)\right]^{3/2}}ds\right)dt\\
&=\frac{1}{k}\sum_{i:c_i^2<1} \int_0^r\int_0^t \frac{\gamma_i^2(1-c_i^2)}{\left[(\gamma_ic_i+s)^2+\gamma_i^2(1-c_i^2)\right]^{3/2}}ds\:dt\\
&=\frac{1}{k}\sum_{i:c_i^2<1} \int_0^r\int_0^t \frac{1-c_i^2}{\gamma_i}\frac{1}{\left[(c_i+\frac{s}{\gamma_i})^2+(1-c_i^2)\right]^{3/2}}ds\:dt.
\end{align*}
Given that 
\[
\left(c_i+\frac{s}{\gamma_i}\right)^2+(1-c_i^2)=\frac{s^2}{\gamma_i^2}+2c_i\frac{s}{\gamma_i}+1\leq \frac{s^2}{\gamma_i^2}+2\frac{s}{\gamma_i}+1=\left(1+\frac{s}{\gamma_i}\right)^2,
\]
we obtain the lower bound
\begin{align*}
F(z) - F(\wh m)&\geq\frac{1}{k}\sum_{i:c_i^2<1} \int_0^r\int_0^t \frac{1-c_i^2}{\gamma_i}\frac{1}{\left[(\frac{s}{\gamma_i}+1)^2\right]^{3/2}}ds\:dt\\
&=\frac{1}{k}\sum_{i:c_i^2<1} \int_0^r\int_0^t \frac{1-c_i^2}{\gamma_i}\frac{1}{(\frac{s}{\gamma_i}+1)^3}ds\:dt\\
&=\frac{1}{k}\sum_{i:c_i^2<1} \int_0^r\int_0^t \frac{\gamma_i^2(1-c_i^2)}{(s+\gamma_i)^3}ds\:dt\\
&=\frac{1}{k}\left(\sum_{j=1}^k \gamma_j^2(1-c_j^2)\right) \int_0^r\int_0^t \sum_{i=1}^k\frac{\gamma_i^2(1-c_i^2)}{\sum_{j=1}^k \gamma_j^2(1-c_j^2)}\frac{1}{(s+\gamma_i)^3}ds\:dt.
\end{align*}
Noting that the inverse cubic function is convex, Jensen's inequality and straightforward integration yields
\begin{align*}
F(z) - F(\wh m) &\geq\frac{1}{k}\left(\sum_{j=1}^k \gamma_j^2(1-c_j^2)\right)\int_0^r\int_0^t \frac{1}{\left(s+\frac{\sum_{i=1}^k \gamma_i^3(1-c_i^2)}{\sum_{j=1}^k \gamma_j^2(1-c_j^2)}\right)^3}ds\:dt\\
&=\frac{1}{2k}\left(\sum_{j=1}^k \gamma_j^2(1-c_j^2)\right) \frac{r^2}{\left(\frac{\sum_{i=1}^k \gamma_i^3(1-c_i^2)}{\sum_{j=1}^k \gamma_j^2(1-c_j^2)}\right)^2\left(r+\frac{\sum_{i=1}^k \gamma_i^3(1-c_i^2)}{\sum_{j=1}^k \gamma_j^2(1-c_j^2)}\right)}.
\end{align*} 
We now observe that
\[
\sum_{i=1}^k \gamma_i^3(1-c_i^2)\leq \sum_{i=1}^k\|\wh m - Y_i\|^3\leq \sum_{i=1}^k \left(\|\wh m \|+\|Y_i\|\right)^3\leq  \sum_{i=1}^k \left(\frac{2}{k} F(0)+\|Y_i\|\right)^3
\]
and also that
\[
\sum_{i=1}^k \gamma_i^2(1-c_i^2) = \sum_{i=1}^k \|\wh m-Y_i\|^2-\left((\wh m-Y_i)^Tu\right)^2=\sum_{i=1}^k\sum_{j=2}^d u_j^T(\wh m-Y_i)(\wh m-Y_i)^T u_j,
\]
where $\{u, u_2, \ldots, u_d\}$ is an orthonormal basis of $\mb{R}^d$. We further notice that
\begin{align*}
\sum_{i=1}^k (\wh m-Y_i)(\wh m-Y_i)^T &= \sum_{i=1}^k (\wh m -\overline{x}+\overline{x}-Y_i)(\wh m-\overline{x}+\overline{x}-Y_i)^T\\
&=k(\wh m-\overline{x})(\wh m-\overline{x})^T +\sum_{i=1}^k(Y_i-\overline{x})(Y_i-\overline{x})^T.
\end{align*}
The Courant-Fischer characterization of eigenvalues gives the inequality
\[
\sum_{i=1}^k \gamma_i^2(1-c_i^2)\geq \sum_{j=2}^d u_j^T\left( \sum_{i=1}^k(Y_i-\overline{x})(Y_i-\overline{x})^T\right)u_j\geq k\sum_{j=2}^d\lambda_j(\widehat{\Sigma}),
\]
where $\{\lambda_j(\widehat{\Sigma})\}_{j=1}^d$ are the eigenvalues of $\wh \Sigma$ listed with multiplicity and in the non-increasing order. We therefore deduce that 
\[
F(z)-F(\wh m)\geq \frac{1}{2} \frac{\sum_{j=2}^d \lambda_j(\widehat{\Sigma}) r^2}{\left(\frac{\frac{1}{k}\sum_{i=1}^k \left(2 m_1+\|Y_i-\overline{x}\|\right)^3}{\sum_{j=2}^d \lambda_j(\widehat{\Sigma})}\right)^2\left(r+\frac{\frac{1}{k}\sum_{i=1}^k \left(2 m_1+\|Y_i-\overline{x}\|\right)^3}{\sum_{j=2}^d \lambda_j(\widehat{\Sigma})}\right)},
\]
thus completing the proof of Theorem \ref{thm:qgc}.

\section{Numerical experiments}
\label{sec:exp}

In this section, we provide experiments to probe the optimality of the dimensional dependence exhibited by the bounds in Theorem \ref{thm:qgc}. We also illustrate the behavior of the geometric median of means estimator on a year's worth of log-returns data from the SNP 500. 

\subsection{Dimension dependence in the bound of Theorem \ref{thm:qgc}}

In this experiment, we generate data by using the real and imaginary parts of random columns of an unnormalized finite Fourier transform matrix. This choice ensures that the data (rows) have norm $\sqrt{d}$. Furthermore, we ensure that the first column is never chosen, so the theoretical mean of the data is $0$, and the covariance matrix of this dataset is always the identity matrix.

\begin{figure}[ht]
    \centering
          \includegraphics[width=0.95\textwidth]{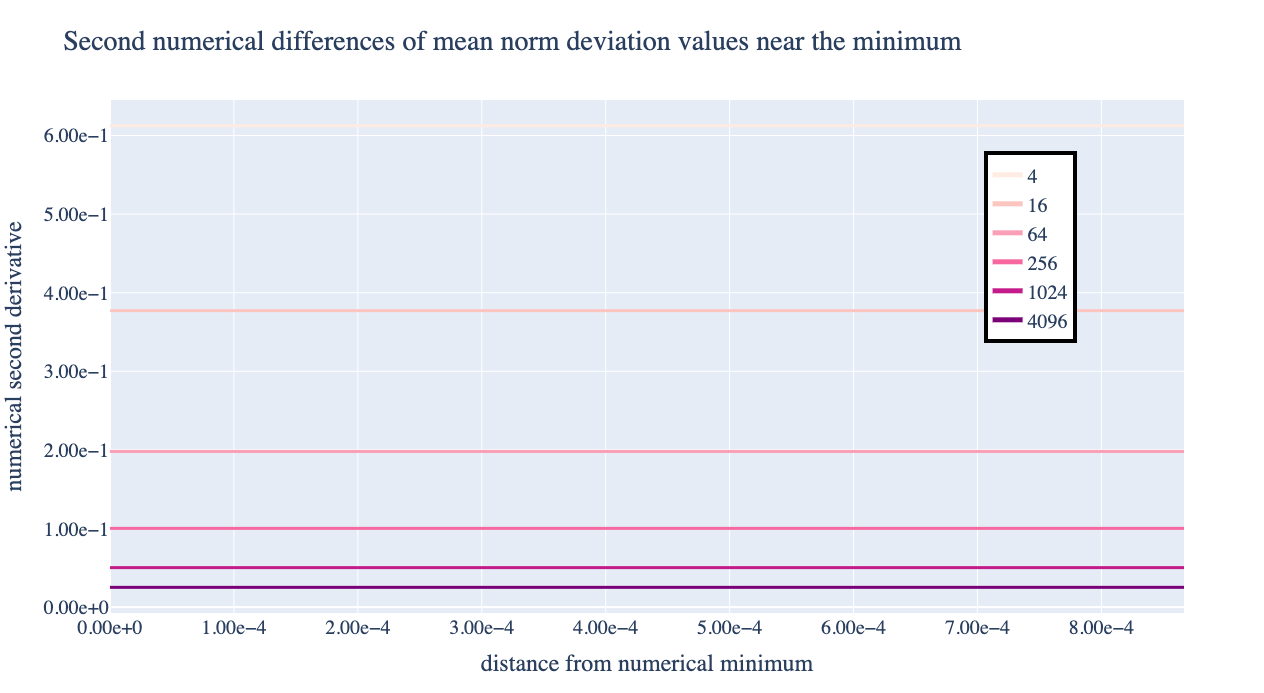}
        \caption{Numerical quadratic growth constants near the geometric median for different dimensions.}
    \label{fig:dim_qgc}
\end{figure}

We consider varying the dimensionality from $d=4, 16, 64, 256, 1024, 4096$, and we always set the number of samples equal to $n=4d$. We compare the numerical mean of the data (the numerical row average) with the theoretical mean of the data ($0$), and we evaluate the mean norm deviation function along the line connecting these two points. In particular, we evaluate at points that are quite close to the numerical mean. To illustrate that the function is essentially a quadratic at these points of evaluation, we observe the numerical second derivatives of these values (Figure \ref{fig:dim_qgc}). This figure shows that these numerical second derivatives appear to be constant, and the values suggest a decay rate on the order of $1/\sqrt{d}$.

\subsection{Extrapolation to future average log-returns}

\begin{figure}[ht]
    \centering
          \includegraphics[width=0.95\textwidth]{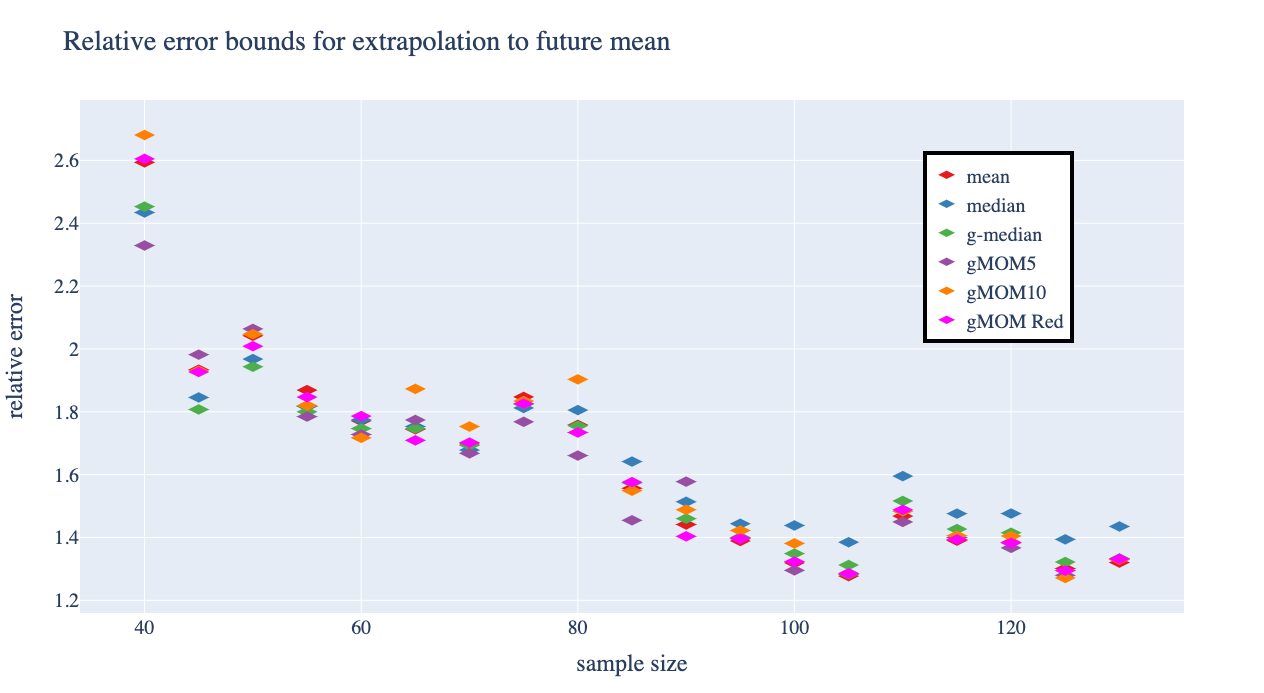}
        \caption{Examples of errors in using various estimators to determine the future average log return.}
    \label{fig:snp}
\end{figure}

Figure \ref{fig:snp} showcases the errors in predicting the future average log return over 100 days using the data accrued for the year so far. The data is from the S\&P 500 index for 2019 calendar year, and the log returns are in terms of the daily opening prices. The dimensionality is the number of symbols observed each day, which is $404$. We test six different procedures: the sample mean (mean), the entry-wise median estimator (median), the geometric median (g-median), the geometric MOM with $k=5$ (gMOM5), the geometric MOM with $k=10$ (gMOM10), and a geometric MOM where each data point is replicated $10$ times and $k=50$ means are used for the geometric median (gMOM Rep). We see that the geometric MOM estimators outperform the other estimators up until there is enough data to make the mean computation more stable, and the geometric MOM estimators outperform the entry-wise median once there is enough data to reveal the bias. 

\bibliographystyle{siamplain}
\bibliography{geomed}

\end{document}

%% file: geomedShared.tex
\usepackage{lipsum}
\usepackage{amsfonts}
\usepackage{graphicx}
\usepackage{epstopdf}
\usepackage{algorithmic}
\ifpdf
  \DeclareGraphicsExtensions{.eps,.pdf,.png,.jpg}
\else
  \DeclareGraphicsExtensions{.eps}
\fi

\usepackage{enumitem}
\setlist[enumerate]{leftmargin=.5in}
\setlist[itemize]{leftmargin=.5in}

\newsiamremark{remark}{Remark}
\newsiamremark{hypothesis}{Hypothesis}
\newsiamremark{example}{Example}
\crefname{hypothesis}{Hypothesis}{Hypotheses}
\newsiamthm{claim}{Claim}

\headers{Geometric Median and Applications}{S. Minsker and N. Strawn}

\title{The Geometric Median and Applications to Robust Mean Estimation\thanks{ Submitted July 2023
\funding{S. Minsker acknowledges support by the National Science Foundation grants DMS CAREER-2045068 and CCF-1908905.}
}}

\author{Stanislav Minsker\thanks{Department of Mathematics, University of Southern California, Los Angeles, CA
  (\email{minsker@usc.edu}).} \and Nate Strawn\thanks{Department of Mathematics and Statistics, Georgetown University, Washington, DC 
  (\email{nate.strawn@georgetown.edu}).}}

\usepackage{amsopn}

\def\R{\mathbb{R}}
